\documentclass[11pt]{article}   % For Latex2e
\usepackage{amssymb,amscd,latexsym}   % For Latex2e
\usepackage[all]{xy}
\usepackage{color}
\LARGE\textwidth=6in
\newcommand{\rar}{\rightarrow}
\newcommand{\lar}{\longrightarrow}
\newcommand{\surjects}{\twoheadrightarrow}

\newtheorem{Theorem}{Theorem}[section]
\newtheorem{Lemma}[Theorem]{Lemma}
\newtheorem{Corollary}[Theorem]{Corollary}
\newtheorem{Proposition}[Theorem]{Proposition}
\newtheorem{Remark}[Theorem]{Remark}
\newtheorem{Example}[Theorem]{Example}

\newtheorem{Definition}[Theorem]{Definition}
\newtheorem{Question}[Theorem]{Question}
\def\sqr#1#2{{\vcenter{\hrule height.#2pt
        \hbox{\vrule width.#2pt height#1pt \kern#1pt
            \vrule width.#2pt}
        \hrule height.#2pt}}}
\def\phi{\varphi}
\def\demo{\noindent{\bf Proof. }}
\def\square{\mathchoice\sqr64\sqr64\sqr{4}3\sqr{3}3}
\def\qed{\hspace*{\fill} $\square$}

\def\xx{{\bf x}}

\def\TT{{\bf T}}
\def\tt{{\bf t}}

\def\uu{{\mathbf{u}}}

\def\fm{{\mathfrak m}}

\def\NN{\mathbb N}

\def\gr#1#2{{\rm gr}\, _{#1}(#2)}
\def\gr{{\rm gr}\,}
\def\hht{{\rm ht}\,}
\def\depth{{\rm depth}\,}
\def\ass{{\rm Ass}\,}

\def\codim{{\rm codim}\,}

\def\ker{{\rm ker}\,}

\def\rk{\rm rank}
\def\hom{\mbox{\rm Hom}}

\def\Rees{{\cal R}}
\def\cl#1{{\cal #1}}

\def\pp{{\mathbb P}}

\begin{document}
%\begin{titlepage}
\begin{center}
{\Large{\bf The Aluffi algebra
}} \footnotetext{
Mathematics   Subject Classification 2000 (MSC2000).
Primary 13A30, 14B05, 14D06, 14H10; Secondary 13B25,
13C15, 13D02, 13F45, 13P10, 14H50
.}\\

\vspace{0.3in}

{\large\sc Abbas  Nasrollah Nejad}\footnote[1]{Supported by a CNPq Graduate Fellowship}
and
{\large\sc Aron  Simis}\footnote[2]{Partially supported by a CNPq Grant.}

\end{center}

%\tableofcontents

\bigskip

\begin{abstract}

We deal with the  {\em quasi-symmetric algebra} introduced
by Paolo Aluffi, here named (embedded) Aluffi algebra.
The algebra is a sort of ``intermediate'' algebra between the symmetric
algebra and the Rees algebra of an ideal, which serves the purpose of
introducing the characteristic cycle of a hypersurface in intersection theory.
The results described in the present paper have an algebraic flavor and
naturally connect with various themes of commutative algebra, such as standard
bases \`a la Hironaka, Artin--Rees like questions, Valabrega--Valla ideals, ideals
of linear type, relation type and analytic spread.
We give estimates for the dimension of the Aluffi algebra and show that,
pretty generally, the latter is equidimensional whenever the base ring is a  hypersurface ring.
There is a converse to this under certain conditions that essentially subsume the setup in Aluffi's theory,
thus suggesting that this algebra will not handle cases other than the singular locus of a hypersurface.
The torsion and the structure of the minimal primes of the algebra are clarified.
In the case of a projective hypersurface the results are more precise and one is naturally led to look at families of
projective plane singular curves to understand how the property of being of linear type deforms/specializes
for the singular locus of a member. It is fairly elementary to show that the singular locus of an
irreducible curve of degree at most $3$ is of linear type.
This is roundly false in degree larger than $4$ and the picture looks pretty wild as we point out
by means of some families.
Degree $4$ is the intriguing case. Here we are able to show that the singular locus of the generic member
of a family of rational quartics, fixing the singularity type, is of linear type.
We conjecture that {\em every} irreducible quartic has singular locus of linear type.

\end{abstract}

\section*{Introduction}

This work is inspired on a  paper of P. Aluffi (\cite{aluffi}) that shows, in the case of a hypersurface,
 how to define a so-called
characteristic cycle in parallel to the well-known conormal cycle in intersection theory.
To accomplish it, Aluffi introduces an intermediate algebra between a symmetric algebra of an ideal
and the corresponding Rees algebra (blowup).

Aluffi has dubbed his construction a {\em quasi-symmetric} algebra.
Since there are many homomorphic images of the symmetric algebra that could equally benefit from this terminology,
we have decided to call it an {\em embedded Aluffi algebra}.
This has the advantage of  indicating that the algebra itself has a more complex
behavior for more general schemes than for
hypersurfaces and, as such, it will often tilt to the other end of the spectrum, namely, become a
honest blowup algebra.

The definition of the algebra is based on taking
ideals $J\subset I\subset R$ in an arbitrary ring, by setting
$${\cal A}_{_{R\surjects R/J}}(I/J):={\cal S}_{R/J}(I/J)
\otimes_{{\cal S}_R(I)}\Rees_R(I).$$

The $R$-embedded Aluffi algebra is functorial in the following sense: let $R\rightarrow R'$ be a ring homomorphism,
let $J'\subset I'\subset R'$ denote the respective images of $J\subset I$ under this map.
If this map induces an isomorphism $R/J\simeq R'/J'$ then it induces a ring surjection
$${\cal A}_{_{R\surjects R/J}}(I/J)\surjects {\cal A}_{_{R'\surjects R'/J'}}(I'/J').$$

Thus, it makes sense to take the inverse limit of such ring surjections.
Letting $A$ denote the common target and ${\mathfrak a}$ the common ideal in the target, Aluffi takes
$${\cal A}({\mathfrak a})\,\colon = \lim_{_{R\surjects A}}\,{\cal A}_{_{R\surjects A}}(\mathfrak a).$$

A point made in his work is that  ${\cal A}({\mathfrak a})$ is actually independent of the choice of the source $R$
(i.e., of the presentation $R/J\simeq A$) provided $R$ is constrained to be regular. Thus, if $R$ is indeed regular then
${\cal A}({\mathfrak a})\simeq {\cal A}_{_{R\surjects A}}(\mathfrak a)$ by the structural map.

Motivated by the case where $R$ is regular, we will study a single member ${\cal A}_{_{R\surjects A}}(\mathfrak a)$
of this inverse system and, accordingly,
omit ``$R$-embedded'' if this causes no confusion.

The overall goal of this work is to study the nature of the algebra {\em ab initio} and then apply it to
a concrete case.
Now, Aluffi focused on the case of a hypersurface -- or, so to say, a Cartier divisor.
Though not explicitly, his work suggests that the application to intersection theory as he had in mind may not
turn out to be suitable for more general varieties.
One of our results explains this insufficiency by showing that if the Aluffi algebra is equidimensional -- actually,
it suffices to know that two of the minimal primes have the same dimension -- then the variety has codimension one.
Moreover, if the Aluffi algebra is actually pure-dimensional (i.e., no embedded primes and equidimensional) then
the variety has to be a Cartier divisor.

Grosso modo the material presented here encloses two sorts of results.
First, one studies the properties of the Aluffi algebra in a quite general ring-theoretic setup,
bringing in some of the typical objects and invariants of commutative algebra.
This will take up the first two sections.
The third section deals with the special case of a projectively embedded hypersurface and its
singular locus (``gradient ideal''), which is the main background in Aluffi's work for this
sort of embedding.

One has a better view of how more structured is the algebra in the case of a
homogeneous equation than in that of its affine companion.
In characteristic zero  the intervenience of the Euler
formula becomes crucial in order to obtain the specifics of the algebra.

A major case, as already pointed out by Aluffi, is the case where the ideal $I$ is of {\em linear type} -- meaning that
${\cl S}_R(I)=\Rees_R(I)$.
As it follows immediately from the definition, this assumption implies that ${\cal A}_{_{R\surjects R/J}}(I/J)={\cal S}_{R/J}(I/J)$.
We show that the converse holds in the case where $J$ is a principal ideal generated by a regular element.

It is our belief that the algebra is relevant on its own and may play a role in situations other than
having $J$ a principal ideal.
Therefore, we deal as much as possible with its structure in the general case -- i.e., when the ideal $J$ has arbitrary codimension.
We find that it is closely related to known themes of commutative algebra, such as
standard bases (\`a la Hironaka),  Artin--Rees number and relation type of an ideal.

One aspect of this to identify the torsion of the Aluffi algebra as the so-called Valabrega--Valla module.
This module -- actually an ideal in the Aluffi algebra -- has been mainly considered in \cite[5.1]{Wolmbook1} in connection
to the situation in which $J$ is a reduction of the ideal $I$.
However, in this case  the structural surjection
$\cl A _{R\surjects A}(I/J)\surjects \Rees_{R/J}(I/J)$ to the relative blowup is actually
an equality in all high degrees, hence the two algebras are finite $R/J$-modules, a case one can dismiss as of
no interest for the present theory, as we are mainly interested in the case where $I$ has a regular element modulo
$J$ -- or at least when $\hht(J)<\hht(I)$ (strict inequality).

An equally meaningful topic is the nature of the associated primes of the Aluffi algebra.
This could throw some light on the summands of the so-called {\em shadow} of the characteristic cycle,
a notion introduced by Aluffi in [loc.cit.] (we thank R. Bedregal for calling our attention
to this matter).
We get pretty close to describing its minimal primes. Since often the algebra is just the symmetric algebra of an ideal,
getting hold of its
associated primes undergoes the same hardship one faces for the latter.
Actually, as we will show, the basic intuition one has about the minimal primes of the symmetric algebra will
work for the Aluffi algebra as well.

Our motivation for the last section comes from Aluffi's quest of the nature of the algebra
in the case that $J$ is generated by the equation of a reduced hypersurface and $I$ defines its singular locus.
Our main interest is to understand its impact in the case of a projective hypersurface and even more modestly,
on the nature of the singularities of plane such curves in low degrees.

Thus, let $J=(f)\subset R$ be a principal ideal,
where $R=k[X_1,\ldots,X_n]$.
We will focus on the Jacobian ideal $I=I_f=(f, \partial f/\partial X_1, \ldots, \partial f/\partial X_n)$.
We are particularly motivated by the problem as to when $I$ is an ideal of linear type.
Now, in general $f$ will not be Eulerian, hence the local number of generators of $I_f$ maybe an early obstruction -- examples of
this sort are easily available.

At the other extreme, if $f$ is Eulerian -- e.g., if $f$ homogeneous in the standard grading of the polynomial ring and and
its degree is not a multiple of the characteristic  -- then it seems like a good bet to expect that $I$ often be of linear type over $R$.
Of course, $I/(f)$ over $R/(f)$ will never be of linear type -- not even generated by analytically independent elements
for that matter -- as the defining equations of the dual variety to $V(f)$ is
a permanent obstruction.

In order to stress the partial derivatives of the homogeneous polynomial $f$ we will call
$I_f$ the {\em gradient} ideal of $f$.
We will assume throughout that char$(k)=0$ or at least that the latter does not divide the degree of $f$.
In this case, by the Euler formula, $f\in I_f$.
One may name the Aluffi embedded algebra in this case the {\em gradient Aluffi algebra} of $f$.

We show that, for a regular element $a\in R$, the Aluffi algebra of a pair
$(a)\subset I\subset R$ is equidimensional. However,  its  structure is still intricate
even when $a$ is the equation of an irreducible projective hypersurface and $I$, its gradient ideal.
A sufficiently tidy case is that of the singular locus of $f$ being set-theoretically a nonempty set of points.
Algebraically,  this translates into  the gradient ideal $I_f$ being a (strict)
almost complete intersection, a situation which is fairly manageable.

We will by and large consider the property of being of linear type for $I_f$ along certain families of plane curves.
Part of the
difficulty of the theory is that, perhaps unexpectedly, the notion of being of linear type is neither kept by
specialization nor by generization.

We discuss the property of being of linear type for the gradient ideal, giving evidence that its behavior
may be rather erratic. The main question is to understand how the nature of the singularities reflect on the algebra
and on its minimal primes, with an eye to the cycle components of Aluffi's characteristic shadow.

Here we are able to show that the singular ideal of the general member
of a family of irreducible rational quartics, fixing the singularity type, is of linear type.
The proof is of some substance as it uses a classification of these curves in terms of quadratic Cremona maps.
We conjecture that any irreducible quartic has gradient ideal of linear type.
In the case of rational quartics, a classification of the possible families allows for a computational verification of
this conjecture. However, we have found no theoretical argument that works for all rational quartics and not just
for the general member of each of these families.

\section{The embedded Aluffi algebra}

\subsection{Preliminaries}

Let $A$ be a ring and ${\mathfrak a}$ an ideal of $A$. Let $R$ be a ring surjecting
onto $A$ and let  $I$ denote the inverse image of ${\mathfrak a}$ in $R$. Note that
the symmetric algebra ${\cal S}_R(I)$ maps surjectively both to the Rees algebra $\Rees_R(I)$
and (by functoriality  of the symmetric algebra) to ${\cal S}_A({\mathfrak a})$.

\begin{Definition}\rm The $R$-{\em embedded
Aluffi algebra} of ${\mathfrak a}$ is defined by
$${\cal A}_{_{R\surjects A}}({\mathfrak a})\colon ={\cal S}_A({\mathfrak a})
\otimes_{{\cal S}_R(I)}
\Rees_R(I).$$
\end{Definition}

We develop a few general preliminaries about the $R$-embedded Aluffi algebra.
The first is a useful presentation that has  already been observed in \cite[Theorem 2.9]{aluffi}
in the context of schemes.

\begin{Lemma}\label{aluffi2rees}
In the above setup, write $I/J:={\mathfrak a}\subset R/J:=A$, where $I\subset R$ is the
inverse image of ${\mathfrak a}$ in $R$.
There are natural $A$-algebra isomorphisms
$${\cal A}_{_{R\surjects A}}({\mathfrak a})\simeq \frac{\Rees_R(I)}{(J,\widetilde{J})\,
\Rees_R(I)}\simeq \bigoplus_{t\geq 0} I^t/JI^{t-1}.$$
where $J$ is in degree $0$ and $\widetilde{J}$ is in degree
$1$. In particular, there is a surjective $A$-algebra homomorphism
${\cal A}_{_{R\surjects A}}({\mathfrak a})\surjects \Rees_A({\mathfrak a})$.
\end{Lemma}
\demo
By the universal property of the symmetric algebra, one sees that
$${\cal S}_{R/J}({I/J})\simeq \frac{{\cal S}_R(I)}{(J,\tilde{J})\,
{\cal S}_R(I)}.$$
Tensoring with $\Rees_R(I)$ gives the first isomorphism.
The second one is now immediate from the definition of $\tilde{J}$.
\qed

\medskip

From the definition and Lemma~\ref{aluffi2rees}, the Aluffi algebra is squeezed as
\begin{equation}\label{squeezing}
{\cal S}_{R/J}(I/J)\surjects {\cal A}_{_{R\surjects R/J}}
(I/J)\surjects \Rees_{R/J}(I/J)
\end{equation}
and is moreover a residue ring of the ambient Rees algebra $\Rees_R(I)$.

If no confusion arises, for a fixed ambient $R$ we will simply refer to this
algebra as the Aluffi algebra of ${\mathfrak a}=I/J$.
Unless stated otherwise, we will assume that $J\subsetneq I\subsetneq R$.

Note that if the ideal $I$ is of linear type -- i.e., the natural surjection
${\cal S}_R(I) \surjects\Rees_R(I)$ is injective -- then trivially
${\cal S}_A(I/J)\simeq {\cal A}_{_{R\surjects A}}(I/J)$.
The following example shows that in general there is no converse to this statement
even when $R$ is a hypersurface domain.
\begin{Example}\rm Let $R=k[x,y,z]=k[X,Y,Z]/(XY-Z^2)$, with $J=(x,z)$ (the ideal of a ruling in the affine cone)
and $I=(x,y,z)$.
Then $R/J\simeq k[Y]$ and $I/J\simeq (Y)$.
Therefore, $I/J$ is of linear type, hence ${\cal S}_{R/J}({I/J})\simeq \cl A_{R\surjects R/J}(I/J)\simeq \Rees_{R/J}(I/J)$.
\end{Example}
This is a particular instance in the following large class: take $(R,\fm)$ to be a non-regular local ring --
or a non-degenerate standard graded algebra
over a field and its irrelevant ideal -- with $J\subset I=\fm$ such that $R/J$ is regular. Then $I/J$ is generated
by a regular sequence on $R/J$, hence is of linear type, while $\fm$ is never of linear type.

It would be interesting to find such examples with $(R,\fm)$ a regular local ring and $J\subset\fm I$.

In the special case where $J$ is a hypersurface, no such examples exist as we now indicate.

\begin{Proposition}\label{lineartype}
Let $R$ be a Noetherian ring and let $I$ denote an ideal. If $a\in I$ is a regular element
then $I$ is of linear type if and only if the natural surjection
$${\cal S}_{R/(a)}({I/(a)})\surjects {\cal A}_{_{R\surjects R/(a)}}(I/(a))$$
is an isomorphism.
\end{Proposition}
\demo
The trivial implication has already been mentioned above.
For the reverse direction, consider an $R$-algebra presentation $S:=R[\TT]\surjects \Rees_R(I)$
based on a set of generators $\mathbf{b}=\{b_1,\ldots, b_n\}$ of $I$.
Write $\mathcal{J}=\bigoplus_{i\geq 1}{\cl J}_i$ for the kernel of this map, where ${\cl J}_i$ stands for
the homogenous part of degree $i$ of $\mathcal{J}$ in the standard grading of $S$.
Note that ${\cl J}_1S\subset \mathcal{J}$  defines likewise the symmetric algebra of $I$ on
$S$, so we need to show that for any $r\geq 0$, $\cl J_r\subset {\cl J}_1S$.

We induct on $r$, the result being trivial if $r=1$. Thus, let $r\geq 2$.
 By Lemma \ref{aluffi2rees} one has
$${\cl J}\subset ({\cl J}_1, \tilde{a}, a)S.$$
Let $F=F(\TT)\in {\cl J}_r$. Then $F=L +\tilde{a}G+aH$
where $L\in {\cl J}_1S_{r-1}$, $G\in S_{r-1}$
and $H\in S_r$.
Note that, if $a=\sum_{j=1}^n c_jb_j$ then $\tilde{a}=\sum_{j=1}^n c_jT_j$, hence $\tilde{a}(\mathbf{b})=a$, i.e.,
evaluating $\tilde{a}$ on the generators of $I$ gives back $a$.
Therefore
\begin{eqnarray*}
0 &=& F(\mathbf{b})=L(\mathbf{b})+\tilde{a}(\mathbf{b})G(\mathbf{b})+a H(\mathbf{b})\\
&=& a\cdot(G+H)(\mathbf{b}),
\end{eqnarray*}
since $L\in \cl J$. As $a$ is a
regular element, $G+H\in{\cl J}$, hence, by homogeneity, $G\in {\cl J}_{r-1}$ and $H\in {\cl J}_r$.

By the inductive hypothesis, $G\in {\cl J}_1S_{r-2}$, hence $\tilde{a}G\in {\cl J}_1S_{r-1}$.
Therefore, $F\in ({\cl J}_1S_{r-1})S+a\,{\cl J}_rS $, thus
showing the equality of ideals ${\cl J}_rS=({\cl J}_1S_{r-1})S+a\,{\cl J}_rS$.
By the graded version of Nakayama's lemma, this
implies that ${\cl J}_rS=({\cl J}_1S_{r-1})S$, as was to be shown.
\qed

\medskip

There is the following consequence, for which we claim no priority.

\begin{Corollary}\label{relative_lineartype}
Let $R$ be a Noetherian ring and let $\{a_1,\ldots,a_m\}\subset R$ be a regular sequence.
If $I\subset R$ is an ideal containing $\{a_1,\ldots,a_m\}$
such that $I/(a_1,\ldots,a_m)$ is
of linear type on $R/(a_1,\ldots,a_m)$ then $I$ is of linear type on $R$.
\end{Corollary}
\demo
Induct on $m$.
For $m=1$, it readily follows from (\ref{squeezing}) and  Proposition~\ref{lineartype}.

Next assume that $m\geq 2$ and write $J=(a_1,\ldots,a_m)$.
Set $\bar{R}:=R/(a_1,\ldots,a_{m-1})$ and, likewise,
$$\bar{J}:=J/(a_1,\ldots,a_{m-1})=(\overline{a_m})\subset \bar{I}:=I/(a_1,\ldots,a_{m-1}).$$
Clearly, $\bar{I}/(\overline{a_m})\simeq I/J$ in $\bar{R}/(\overline{a_m})\simeq R/J$.
Therefore, the assumption that $I/J$ is of linear type on $R/J$ implies that $\bar{I}/(\overline{a_m})$
is of linear type on $\bar{R}/(\overline{a_m})$,  where $\overline{a_m}$
is a regular element in $\bar{R}$.
By the first part, $\bar{I}$ is of linear type on $\bar{R}$.
Then, by the inductive hypothesis, $I$ is of linear type on $R$.
\qed

\subsection{Dimension}

A few routine statements follow from the preliminaries of the previous subsection.

\begin{Proposition}\label{dimension_bounds}
Let $J\subsetneq I\subsetneq R$ be ideals of the Noetherian ring $R$.
\begin{enumerate}
\item[{\rm (a)}] If $J$ has a regular element then
$\dim \cl A _{R\surjects R/J}(I/J)\leq \dim R$.
\item[{\rm (b)}] If $I/J$ has a regular element then
$$\min\{\dim R+1\,,\, \dim {\cal S}_{R/J}(I/J)\,\}\geq \dim \cl A _{R\surjects
R/J}(I/J)\geq \dim R/J+1.$$
\end{enumerate}
\end{Proposition}
\demo (a) Since $\Rees_R(I)$ is $R$-torsionfree, one has $\hht J\Rees_R(I)\geq 1$.
Therefore
$$\dim \cl A _{R\surjects R/J}(I/J)\leq \dim \Rees_R(I)/J\Rees_R(I)\leq \dim R+1-1
=\dim R.$$

(b) This follows immediately from (\ref{squeezing}) by the well-known
dimension formula for the Rees algebra of an ideal containing a regular element.
\qed

\begin{Remark}\rm In (a) this is all one can assert in such generality because if, e.g., a power
of the ideal $I$ is
contained in $J$, then $\dim \cl A _{R\surjects R/J}(I/J)= \dim R/J$.
\end{Remark}

Perhaps less routine is the following result.

\begin{Theorem}\label{hypersurface}
Let $R$ be a catenary, equidimensional and equicodimensional Noetherian ring
and let $I\subsetneq R$ be an ideal containing a regular element $a$.
Then $\cl A _{R\surjects R/J}(I/(a))$ is equidimensional and
$\dim \cl A _{R\surjects R/J}(I/(a))= \dim R$.
\end{Theorem}
\demo
Under the assumptions on $R$, $I/(a)$ and $a$, one can apply Proposition~\ref{dimension_bounds}, (a),
and the right hand inequality of (b) to conclude that $\dim \cl A _{R\surjects R/J}(I/(a))=\dim R-1+1=\dim R$.

To prove the equidimensionality part,
we will show that $\cl A _{R\surjects R/(a)}(I/(a))$
is equidimensional locally at every prime ideal  ${\cl P}\subset \Rees_R(I)$ in its support.
Localizing first at ${\cl P}\cap R$ in the base ring one can assume that $(R,\fm)$ is local, with ${\cl P}\cap R=\fm$
and $I\subset \fm$.
Now, ${\cl M}=(\fm,Iu)\subset \Rees_R(I)\subset R[Iu]$ is not a minimal prime of $\cl A _{R\surjects R/(a)}(I/(a))$. This is because
the Aluffi algebra is graded, with grading induced from $\Rees_R(I)$, hence $\cl M$ would actually be its unique
associated prime, which is impossible as $\dim R\geq 1$.
Thus, for the purpose of showing equidimensionality, we may assume that ${\cl P}$ is a homogeneous ideal properly contained in
${\cl M}$.

Let $I=(b_1,\ldots,b_n)$.
Note that in the present situation, one has by Lemma~\ref{aluffi2rees}:
 $$\cl A _{R\surjects R/(a)}(I/(a))\simeq \Rees_R(I)/(a,\tilde{a}).$$
 Write $\Rees_R(I)=R[b_1u,\ldots, b_nu]\subset R[u]$, so that $\tilde{a}=\sum_{j=1}^n c_jb_ju$,
 for suitable $c_j\in R$.

 Suppose first that $(Iu)\not\subset {\cl P}$.
 Say, $b_1u\notin {\cl P}$. Localizing at ${\cl P}$ yields

\begin{eqnarray*}\cl A _{R\surjects R/(a)}(I/(a))_{\cl P}&\simeq &R[Iu]_{\cl P}/(a,\tilde{a})_{\cl P} \simeq
R\left[\frac{I}{b_1}\,, b_1u, (b_1u)^{-1}\right]_{{\cl P}'}\biggl/\left(a,c_1+\sum_{j=2}^m c_j\frac{b_j}{b_1}\right)_{{\cl P}'}\\
&=& R\left[\frac{I}{b_1}\,, b_1u, (b_1u)^{-1}\right]_{{\cl P}'}\biggl/\left(a,\frac{a}{b_1}\right)_{{\cl P}'}\\[8pt]
&=& R\left[\frac{I}{b_1}\,, b_1u, (b_1u)^{-1}\right]_{{\cl P}'}\biggl/\left(\frac{a}{b_1}\right)_{{\cl P}'},
\end{eqnarray*}
where ${\cl P}'$ denotes the corresponding image of ${\cl P}$.
The rightmost ring above is a factor ring of a catenary, equidimensional and equicodimensional
ring by a principal ideal generated by a regular element, hence it is equidimensional and so is
$\cl A _{R\surjects R/(a)}(I/(a))_{\cl P}$ too.

Suppose now that $(Iu)\subset {\cl P}$. Then $\fm\not\subset {\cl P}$ since ${\cl P}\subsetneq {\cl M}$,
hence $p:={\cl P}\cap R\subsetneq \fm$.
Note that $p$ is a prime containing $a$.

If $I\not\subset p$ then $\cl A _{R\surjects R/(a)}(I/(a))_{\cl P}$ is a localization of the ring
$$R_p[I_pu]/(a,\tilde{a})=R_p[u]/(a,au)=R_p[u]/(a)$$
and we conclude as above.
If $I\subset p$ then $\cl A _{R\surjects R/(a)}(I/(a))_{\cl P}$ is a localization of the Aluffi algebra
$\cl A _{R_p\surjects R_p/(a)}(I/(a))$ and we conclude by induction on $\dim R$.

\qed

\medskip

A geometric version of Theorem~\ref{hypersurface} case is stated in \cite[Corollary 2.18]{aluffi}.

We will have more to say about the equidimensionality of the Aluffi algebra in subsequent sections.

\subsection{Local or graded case}

In this part we assume that $(R,\fm)$ is a Noetherian local ring and its maximal ideal or a standard
graded algebra over a field and its maximal irrelevant ideal.
Throughout $R/\fm$ is an infinite field.

Let $J\subset I\subset \fm$.
We confront ourselves with two quite opposite situations, namely, when $J\subset \fm I$ and when $J$ contains
minimal generators of $I$.
Note that if $J$ is a reduction of $I$ then $J\subset \fm I$
would entail $I^t=JI^{t-1}\subset \fm I^t,$ for $t>>0$,
hence $I^t=\{0\}$, i.e., $I$ would be nilpotent.

Drawing upon a terminology of geometry, let us agree to say that the pair $J\subset I$ of ideals is
{\em non-degenerate} if $J\subset \fm I$. If on the other extreme, $J\subset I$ is generated by a subset of
minimal generators of $I$, we may call the pair $J\subset I$ {\em totally degenerate}.
The latter case can usually be disposed of by a standard argument (see Proposition~\ref{residual_gens}).

\medskip

We recall that the {\em analytic spread} of $I$, denoted $\ell(I)$, is the dimension of
the $R/\fm$-standard algebra $\Rees_R(I)/\fm \Rees_R(I)$.
It can be shown that $\ell(I)$ coincides with the number of minimal generators of a least
possible reduction of $I$, but we shall have no occasion to use this result.
The behavior of $\ell(I)$ in face of other numerical invariants related to $I$ is as follows:
$$\hht I\leq \ell(I)\leq \min\{\mu(I),\dim R\},$$
where $\mu(I)$ denotes minimal number of generators.
We will say that $I$ has {\em maximal analytic spread} if $\ell(I)=\dim R$. Note that this forces
$\mu(I)\geq \dim R$.

\begin{Proposition}\label{lying_deep}
Let $(R,\fm)$ be as above with $R/\fm$ infinite.
Suppose that $J\subset I$ is a non-degenerate pair and that $J$ has  a regular element.
Then
\begin{enumerate}
\item[{\rm (i)}] $\ell(I)\leq \dim \cl A _{R\surjects R/J}(I/J)\leq\dim R${\rm ;}
in particular, if $I$ has maximal analytic spread then
$\dim \cl A _{R\surjects R/J}(I/J)=\dim R$.
\item[{\rm (ii)}] If $I$ has maximal analytic spread and, moreover, $\mu(I)=\dim R$,
then $\fm\Rees_R(I)$ is a minimal prime of $\cl A _{R\surjects R/J}(I/J)$
of maximal dimension.
\item[{\rm (iii)}]  If $J\subset I^2$ then
$\dim \cl A _{R\surjects R/J}(I/J)=\dim R$.
\end{enumerate}
\end{Proposition}
\demo (i)
$J\subset \fm I$ implies $JI^{t-1}\subset \fm I^t$ for every $t\geq 0$.
 This yields a surjective homomorphism
 $\cl A_{R\surjects R/J}(I/J)\surjects \Rees_R(I)/\fm \Rees_R(I)$, from which follows the leftmost
 inequality.

 The other inequality stems from Proposition~\ref{dimension_bounds}, (a).

 \medskip

 (ii) The assumption $\ell(I)=\mu(I)=\dim R$ implies that $I$ is generated by
 analytically independent elements and
 the latter entails that $\Rees_R(I)/\fm \Rees_R(I)$ is a polynomial ring over $R/\fm$.
 In particular, $\fm\Rees_R(I)$ is a prime ideal of $\Rees_R(I)$.
Since $JI^{t-1}\subset \fm I^t\subset \fm$ then $(J,\tilde{J})\subset \fm\Rees_R(I)$
as ideals of $\Rees_R(I)$.
Therefore  $\cl A _{R\surjects R/J}(I/J)/\fm\, \cl A _{R\surjects R/J}(I/J)\simeq \Rees_R(I)/\fm \Rees_R(I)$, hence $\fm\cl A _{R\surjects R/J}(I/J)$ is a prime ideal with
$$\dim\cl A _{R\surjects R/J}(I/J)/\fm \,\cl A _{R\surjects R/J}(I/J)=\dim R=\dim\cl A _{R\surjects R/J}(I/J)$$
by the first part.

\medskip

(iii) Write $\gr_I(R)$ for the associated graded ring of $I$.
Since $J\subset I^2$, one has $JI^{t-1}\subset I^{t+1}$ for every $t\geq 0$.
This yields a surjective homomorphism
 $\cl A_{R\surjects R/J}(I/J)\surjects \gr_I(R)$, showing that
 $\dim \cl A _{R\surjects R/J}(I/J)\geq \dim R$.
 The reverse inequality follows from Proposition~\ref{dimension_bounds}, (a).
\qed

\medskip

We wrap up with a comment on the last result. Namely, we actually have
$$\dim \cl A _{R\surjects R/J}(I/J)\geq \max\{\ell(I), \dim R/J+1\}$$
provided $I$ has a regular element modulo $J$.
The interesting case is when $\ell(I)\geq \dim R/J+1$. If, say, $R$ is catenary and
equidimensional, it would entail
$$\dim \cl A _{R\surjects R/J}(I/J)
\geq \frac{\dim R+1}{2}.$$

\section{Structural properties}

In this section one looks more closely at the internal structure of the Aluffi algebra
and relate some of the elements of this structure to well-known notions in ideal theory.

\subsection{Torsion and minimal primes}

By Lemma~\ref{aluffi2rees} one has
$$\cl A _{R\surjects R/J}(I/J)\simeq
\Rees_R(I)/(J,\tilde{J})\Rees_R(I)=\bigoplus_{t\geq 0} I^t/JI^{t-1}.$$
Since $\Rees_{R/J}(I/J)=\bigoplus_{t\geq 0} (I^t,J)/J\simeq \bigoplus_{t\geq 0} I^t/J\cap I^{t}$,
 the kernel of  the natural surjection ${\cal A}_{_{R\surjects R/J}}
(I/J)\surjects \Rees_{R/J}(I/J)$ is the homogeneous ideal
\begin{equation}\label{vava}
{\cal V}\kern-5pt {\cal V}_{J\subset I}:=\bigoplus_{t\geq 2} \frac{J\cap
I^t}{JI^{t-1}},
\end{equation}
dubbed as the {\em module of Valabrega--Valla} (see
\cite{VaVa}, also \cite[5.1]{Wolmbook1})

\medskip

We retrieve a result of Valla (\cite[Theorem 2.8]{Valla2}):

\begin{Corollary}\label{valla_lineartype}
Let $J\subset I\subsetneq R$ be ideals of the local ring $R$.
If $I/J$ is of linear type over $R/J$ {\rm (}e.g., if $I$ is generated by
a regular sequence modulo $J${\rm )} then $J\cap I^t=JI^{t-1}$
for every positive integer $t$.
\end{Corollary}
\demo This follows immediately from the structural ``squeezing'' (\ref{squeezing}).
\qed

\medskip

 Note that the assumption in
\cite[Proposition 3.10]{Trento} to the effect that $I$ be of linear type over $R$
does not intervene in the above statement.

Here is a useful explicit situation, where we write $I=(J,\mathfrak{a})$, with no particular care
for minimal generation.

\begin{Proposition}\label{residual_gens}
Let $I=(J,\mathfrak{a})$.
If $J\cap \mathfrak{a}^t\subset J\mathfrak{a}^{t-1}$, for every $t\geq 0$ then
$\cl A _{R\surjects R/J}(I/J)\surjects \Rees_{R/J}(I/J)$ is an isomorphism.
\end{Proposition}
\demo One has:
$$J\cap I^t=J\cap (J,\mathfrak{a})^t=J\cap (J(J,\mathfrak{a})^{t-1}, \mathfrak{a}^t)=
J(J,\mathfrak{a})^{t-1}+J\cap \mathfrak{a}^t\subset
JI^{t-1}+J\mathfrak{a}^{t-1}\subset JI^{t-1}.$$
\qed

\begin{Remark}\rm In the notation of the previous proposition, one of the main results of
\cite{Huneke} is that if $\mathfrak{a}$ is generated by a $d$-sequence modulo $J$ then
the assumption of the proposition is fulfilled.
Therefore, under the hypothesis of [{\em loc.cit.}],
the surjection $\cl A _{R\surjects R/J}(I/J)\surjects \Rees_{R/J}(I/J)$ is an isomorphism.
This result, however, is a special case of Corollary\ref{valla_lineartype} if one uses that
an ideal generated by a $d$-sequence is of linear type.
Of course, the proof of this fact requires some non-trivial work on itself and is previous
to the later results, such as \cite{Trento}.
\end{Remark}

When $J=(a)$ is a principal ideal, one has
a result somewhat subsumed in the spirit of \cite{Valla2}.

\begin{Proposition} Let $\mathfrak{a}$ be an ideal
 in the ring $R$ and let $a\in R$ be an
element such that $\mathfrak{a}^t:a=\mathfrak{a}^t$ for every integer $t\geq 0$. Then the inclusion
$(a)\subset (a,\mathfrak{a})$ induces an isomorphism
${\cal A}_{_{R\surjects R/(a)}}
((a,\mathfrak{a})/(a))\simeq \Rees_{R/(a)}((a,\mathfrak{a})/(a))$.
\end{Proposition}
\demo
The assumption means that $(a)\cap \mathfrak{a}^t=a\mathfrak{a}^t$ for every $t\geq 0$, hence
 $(a)\cap (a,\mathfrak{a})^t=a(a,\mathfrak{a})^{t-1}+(a)\cap \mathfrak{a}^t=a(a,\mathfrak{a})^{t-1}$ for $t>0$.
 \qed

\medskip

The Valabrega--Valla module  gives the torsion in as many cases as the ones in which
the Rees algebra is the symmetric algebra modulo torsion.

\begin{Proposition}\label{torsion}
Let $J\subset I\subsetneq R$ be ideals of the Noetherian ring $R$.
If $I/J$ has a regular element then the $R/J$-torsion of the embedded Aluffi algebra
of $I/J$ is the kernel of the natural surjection ${\cal A}_{_{R\surjects R/J}}
(I/J)\surjects \Rees_{R/J}(I/J)$.
\end{Proposition}
\demo Consider the general elementary observation: given a ring $A$ and $A$-modules
$$N\surjects M\surjects K$$
such that the $A$-torsion of $N$ is the kernel of the composite
$N\surjects K$ then the $A$-torsion of $M$ is the kernel of
$M\surjects K$.
We apply this to the situation in (\ref{squeezing}), by recalling that if ${\mathfrak a}\subset A$ has
a regular element in the ring $A$ then the $A$-torsion of the symmetric algebra $\cl S _A({\mathfrak a})$
is the kernel of the natural surjection $\cl S _A({\mathfrak a})\surjects \Rees_A({\mathfrak a})$.
\qed

\medskip

Recall that, given a ring $S$, an ideal $\mathfrak{b}\subset S$ and an $S$-module $E$, one
denotes by H$^0_{\mathfrak{b}}(E)$ the zeroth local cohomology of $E$ with respect to
$\mathfrak{b}$. One has
$$ {\rm H}^0_{\mathfrak{b}}(E)\simeq E:\mathfrak{b}^{\infty}:=
\{\,\epsilon\in E\,|\, \exists n\geq 0,\, \mathfrak{b}^n\epsilon=\{0\}\,\}.$$

\begin{Corollary}\label{torsion_and_H0}
Let $J\subset I\subsetneq R$ be ideals of the Noetherian ring $R$.
If $I/J$ has a regular element then ${\cal V}\kern-5pt {\cal V}_{J\subset I}={\rm H}^0_{I/J}(\cl A)$,
where $\cl A$ denotes the Aluffi algebra and ${\cal V}\kern-5pt {\cal V}_{J\subset I}$ is
the Valabrega--Valla module as introduced above.
\end{Corollary}
\demo By Proposition~\ref{torsion}, ${\cal V}\kern-5pt {\cal V}_{J\subset I}$ is the $R/J$-torsion of
$\cl A$. On the other hand, localizing at primes of the base $R/J$ not containing $I/J$ makes the
surjection $\cl S _{R/J}(I/J)\surjects\Rees_{R/J}(I/J)$ an isomorphism, hence also the surjection
$\cl A\surjects \Rees_{R/J}(I/J)$. Therefore, $\cl A$ is torsionfree locally at those primes.
Since $I/J$ has regular elements, the result follows easily (see, e.g., \cite[Lemma 5.2]{elam99}).
\qed

\begin{Remark}\label{exponent}\rm The last result says, in particular, that there exists an integer $k\geq 0$ such that
$I^k\,(J\cap I^t)\subset JI^{t-1}$ for every $t\geq 1$.
Later on we will relate such an exponent to the so-called Artin--Rees number.
\end{Remark}

By the same principle, one can get a hold of the minimal primes of the Aluffi algebra.
Quite generally, to any ideal ${\mathfrak a}\subset R$ we associate its extended--contracted ideal
$$\tilde{{\mathfrak a}}:={\mathfrak a}R[u]\cap R[Iu]=\sum_{t\geq 0} ({\mathfrak a}\cap I^t)u^t$$
in the Rees algebra $\Rees_R(I)\simeq R[Iu]\subset R[u]$ ($u$ a variable over $R$).

\begin{Proposition}\label{minimal_primes}
Let $J\subset I\subsetneq R$ be ideals of the Noetherian ring $R$.
Any minimal prime $\wp$ of $\cl A _{R\surjects R/J}(I/J)$ on $\Rees_R (I)$  is either of the form $\wp=\tilde{p}$ for
some minimal prime of $R/J$ on $R$, or else has the form $(q,\wp_+)$ where $q:=\wp\cap R$ contains a minimal
prime of $R/I$ on $R$ and $\wp_+=\wp\cap (Iu)$.
\end{Proposition}
\demo By Corollary~\ref{torsion_and_H0} -- rather by its proof -- a power of $I$ annihilates the kernel of
$\cl A _{R\surjects R/J}(I/J)\surjects \Rees_{R/J}(I/J)$ lifted to $\Rees_R(I)$ -- call it $\cl K$.
If $\wp\subset \Rees_R(I)$ is a minimal prime of $\cl A _{R\surjects R/J}(I/J)$ it follows that $\wp$ contains either
$\cl K$ or $I$. In the first case, it contains a minimal prime of $\Rees_R(I)/{\cl K}\simeq \Rees_{R/J}(I/J)$ hence
must be a minimal prime of the latter on $\Rees_R(I)$.
But, it is well known that the above extending-contracting operation induces a bijection between the minimal primes of $R/J$
on $R$ and the minimal primes of $\Rees_{R/J}(I/J)$ on $\Rees_R(I)$.

In the case $I\subset \wp$, since $\wp$ is homogeneous in the natural $\NN$-grading of $\Rees_R(I)$, then it is clear
that $\wp=(q,\wp_+)$, where $\wp_+=\wp\cap \Rees_R(I)_+$; obviously, $q$ contains a minimal prime of $R/I$ on $R$.
\qed

\begin{Remark}\rm
Note that if $\wp\subset \Rees_R(I)$ is a minimal prime of $\cl A _{R\surjects R/J}(I/J)$ containing $I$
then $\wp_+$ behaves erratically: it can actually be zero in certain cases (see, e.g., Proposition~\ref{lying_deep}, (ii)).
On the other hand, its contraction $q\subset R$ may turn out to be an embedded associated prime of $R/I$ and not
a minimal one (see Section~\ref{gradients}).
\end{Remark}

\begin{Corollary}\label{equidimension}
Let $J\subset I\subsetneq R$ be ideals of a Noetherian ring $R$ such that
$I/J$ has a regular element.
If $\cl A _{R\surjects R/J}(I/J)$ is equidimensional then
$$\dim \cl A _{R\surjects R/J}(I/J)=\dim R/J+1.$$
\end{Corollary}
\demo This readily follows from Proposition~\ref{minimal_primes} and the
known value of the dimension of $\Rees_{R/J}(I/J)$ under the present assumption on $I/J$.
\qed

\medskip

Equidimensionality of $\cl A _{R\surjects R/J}(I/J)$ if $J$ has codimension $\geq 2$ may be quite rare.
The next result shows that, at least in the local or graded case, pure-dimensionality is really infrequent
under this assumption.

\begin{Proposition}\label{pure_dimension}
Let $(R,\fm)$ be a Noetherian local ring and its maximal ideal $\fm$ or a standard
graded algebra over a field
and its maximal irrelevant ideal $\fm$. Assume that $R/\fm$ infinite.
Let $J\subset I\subset \fm$, with $J$ having a regular element and $I$ having a regular element
modulo $J$.
Suppose that:
\begin{itemize}
\item[{\rm (i)}] $J\subset I\subset \fm$ is a non-degenerate pair and $\ell(I)=\mu(I)=\dim R\,${\rm ;} or else
\item[{\rm (ii)}] $J\subset I^2$.
\end{itemize}
If $\cl A _{R\surjects R/J}(I/J)$ is equidimensional then $J$ has height one.
If $\cl A _{R\surjects R/J}(I/J)$ is pure-dimensional then
$J$ is an ideal of pure height one {\rm (}hence, principal if $R$ is regular{\rm)}.
\end{Proposition}
\demo
By  Proposition~\ref{minimal_primes}, $\Rees _{R/J}(I/J)$ is equidimensional. In particular, one has
$$\dim \cl A _{R\surjects R/J}(I/J)=\dim \Rees _{R/J}(I/J)=\dim R/J+1\leq\dim R-\hht J+1.$$
On the other hand, by Proposition~\ref{lying_deep}, (ii) or (iii), $\dim \cl A _{R\surjects R/J}(I/J)=\dim R$.
Therefore, $\hht J\leq 1$ (hence $\hht J=1$ since $J$ has a regular element by hypothesis).

Now, assume that $\cl A _{R\surjects R/J}(I/J)$ is pure-dimensional.
Let $p\in\ass_R(R/J)$ have height $\geq 2$. By  Proposition~\ref{minimal_primes}
and the the pure-dimensionality of $\cl A _{R\surjects R/J}(I/J)$, the prime
$\tilde{p}$ satisfies
\begin{eqnarray*}\dim \cl A _{R\surjects R/J}(I/J)&=& \dim \cl A _{R\surjects R/J}(I/J)/\tilde{p}
=\dim \Rees_R(I)/\tilde{p}=\dim \Rees_{R/p}((p,I)/p)\\
&=&\dim R/p+1\leq \dim R-\hht p+1
\leq \dim R-2+1=\dim R-1.
\end{eqnarray*}
Again Proposition~\ref{lying_deep}, (ii) or (iii)  gives a contradiction.
\qed

\bigskip

 More difficult is to get hold of non-trivial embedded primes of $\cl A _{R\surjects A}(I/J)$.
In the case where $J$ is the ideal of a homogeneous hypersurface in projective space there are often embedded primes
containing the irrelevant ideal.

We wrap up with the following
\begin{Question}\rm
Suppose as above that $J=(a)\subset \fm I$, $I/(a)$ has a regular element and $I$ has maximal analytic spread.
To what extent can we assert that, conversely, $\cl A _{R\surjects A}(I/J)$  is pure-dimensional?

This seems to be the case in a variety of situations such as the one considered in Section~\ref{gradients}.
\end{Question}

\subsection{Relation to the Artin--Rees number}

A close associate to ${\cal V}\kern-5pt {\cal V}_{J\subset I}$ is the well-known ideal
$$\ker(\gr_I(R)\surjects \gr_{I/J}(R/J) )=\bigoplus_{t\geq 0} ( I^{t+1}+J\cap I^t)/ I^{t+1},$$
generated by the $I$-initial forms of elements of $J$. Recall that the
$I$-initial form of an element $a\in R$ is the residue class $a^*$ of $a$ in
$I^{\nu(a)}/I^{\nu(a)+1}$ where $\nu(a)$ is the $I$th order of $a$ (i.e.,
$a\in I^{\nu(a)}\setminus I^{\nu(a)+1}$, setting $\nu(a)=\infty$ and $a^*=0$ if $a\in \cap_{t\geq 0}I^t$).

A set of elements of $J$ is called an $I$-{\em standard base} of $J$ if
their initial forms in $\gr_I(R)$ generate the above ideal. If $R$ is Noetherian local then an $I$-standard base of $J$ is a generating set
of $J$ (see \cite[Lemma 6]{Hironaka}). We will shorten $\nu(a_i)$ to $\nu_i$ if $a_i$
is sufficiently clear from the context.

The following basic result will be used throughout.

\begin{Theorem}{\rm (\cite{VaVa})}\label{standardbase}
Let $J=(a_1,\ldots,a_m)$ be an ideal of ring $R$. Then
$\{a_1,\ldots,a_m\}$ is an $I$-standard base of $J$ if and only if
$$J\cap I^t=\sum_{i=1}^{m}a_iI^{t-\nu_i}$$
for every positive integer $t$.
\end{Theorem}

\medskip

This result implies immediately:

\begin{Proposition}\label{inclusion}
Let $R$ be a Noetherian local ring and $J\subset I$ be ideals of $R$. Let
$\{a_1,\ldots,a_m\}$ be an $I$-standard base of $J$ such that
$1\leq\nu_1\leq\nu_2\leq\dots\leq\nu_m$.
\begin{enumerate}
\item[{\rm (a)}] The surjection ${\cal A}_{_{R\surjects R/J}}
(I/J)\surjects \Rees_{R/J}(I/J)$ is an isomorphism if and only if $\nu_m=1$ {\rm (}i.e., $\nu_1=\cdots =\nu_m=1${\rm )}
\item[{\rm (b)}] More generally
$$\mathfrak{R} \cap I^{\nu_m-1}\Rees_R(I)\,\subset \,\mathfrak{A}\,\subset\, \mathfrak{R} \cap I^{\nu_1-1}\Rees_R(I),$$
where $\mathfrak{R} =\bigoplus_{t\geq 0}J\cap I^t$ and $\mathfrak{A}=\bigoplus_{t\geq 0}JI^{t-1}$.
\item[{\rm (c)}]   If $\nu_1>1$ and $J$ is not contained in any minimal prime of $R$, then
$\dim \cl A _{R\surjects R/J}(I/J)=\dim R$.
\end{enumerate}
\end{Proposition}
\demo (a) One direction is obvious from Theorem~\ref{standardbase}.
For the reverse implication, assume that $J\cap I^t=JI^{t-1}$ for
every positive integer $t$. By the above remark, one may assume that  $a_1,\dots,a_m$ is a minimal set of
generators of $J$. If for some $i$, $a_i\in I^2$ then by hypothesis $a_i\in
JI$, which clearly contradicts the minimality of $a_1,\dots,a_m$.

(b) By definition, we want to show the two inclusions
$$J\cap I^{t+\nu_m-1}\subset JI^{t-1}\subset J\cap I^{t+\nu_1-1}$$
as subideals of $J\cap I^t$, for every $t\geq 1$.

This is however a straightforward consequence of Theorem~\ref{standardbase} as one has thereof
$$J\cap I^{\nu_m+t-1}=\sum_{i=1}^{m}a_iI^{\nu_m+t-1-\nu_i}\subset
\sum_{i=1}^{m}a_iI^{t-1}=JI^{t-1},$$
and similarly
$$J\cap I^{\nu_1+t-1}=\sum_{i=1}^{m}a_iI^{\nu_1+t-1-\nu_i}\supset
\sum_{i=1}^{m}a_iI^{t-1}=JI^{t-1}.$$

(c) Quite generally, for a positive integer $r$ one has
\begin{eqnarray}
\nonumber \dim \Rees_R(I)/\mathfrak{R} \cap I^r\Rees_R(I)&\geq &\max\{\dim \Rees_R(I)/\mathfrak{R} \ ,\ \dim
\Rees_R(I)/I^r\Rees_R(I)\}\\
\nonumber &=& \max\{\dim \Rees_{R/J}(I/J),\ \dim \gr_I(R)
   \}=\dim R.
\end{eqnarray}
On the other hand,  since $J$ is assumed to be of positive height,
by Proposition~\ref{dimension_bounds} (a) one has $\dim \cl A _{R\surjects R/J}(I/J)\leq \dim R$.
The result follows now at once.
\qed

\bigskip

We close with yet another condition for
the surjection $\cl A _{R\surjects R/J}(I/J)\surjects \Rees_{R/J}(I/J)$ to be an isomorphism.

For this recall that, pretty generally,  given ideals $J,I\subset R$ the {\em Artin--Rees number}
of $J$ relative to $I$ is the integer
$$\min\{k\geq 0\,|\, J\cap I^t=(J\cap I^k)I^{t-k}\;\forall\; t\geq k\}.
$$

We observe that if $J\subset I\subsetneq R$, where $R$ is Noetherian and $J$ has regular elements
then the Artin--Rees number of $J$ relative to $I$ is $\geq 1$.

\begin{Proposition}\label{artinrees_vava}
Let $J\subset I\subsetneq R$ be ideals of a Noetherian ring $R$
and let $k\geq 1$ be an upper bound for the Artin--Rees number of $J$ relative to $I$, i.e., $J\cap I^t=(J\cap I^k)I^{t-k}\;\forall\; t\geq k$.

Then $I^{k-1}$ annihilates the kernel of the surjection$\cl A _{R\surjects R/J}(I/J)\surjects \Rees_{R/J}(I/J)$.
Moreover, the  latter is an isomorphism if and only if the Artin--Rees number of
   $J$ relative to $I$ is $1$.
\end{Proposition}
\demo
One has $(J\cap I^t)I^{k-1}=(J\cap I^k)I^{t-k}I^{k-1}= (J\cap I^k)I^{t-1}\subset JI^{t-1}$ for $t\geq k$.
On the other hand, for $t\leq k-1$, one has $I^{k-1}\subset I^{t-1}$, hence
$(J\cap I^t)I^{k-1}\subset (J\cap I^t)I^{t-1}\subset JI^{t-1}$.

The second assertion is clear.
\qed

\medskip

More generally:

\begin{Lemma}\label{artinrees_isomorphism}
Let $J\subset I\subset R$ be ideals of a ring $R$.
Assume that $\ell$ is an upper bound for the Artin--Rees number of $J$ relative to $I$
such that $J\cap I^{t}=JI^{t-1}$ for every $t\leq \ell$.
Then $\cl A _{R\surjects R/J}(I/J)\surjects \Rees_{R/J}(I/J)$ is an isomorphism.
\end{Lemma}
\demo By assumption $J\cap I^t=I^{t-\ell}(J\cap I^{\ell})$ for $t\geq \ell$. Now use
the assumed equality $J\cap I^{\ell}=JI^{\ell-1}$ to get $J\cap I^t=JI^{t-1}$ for $t\geq \ell$.
\qed

\medskip

Given a ring $A$ and an ideal $\mathfrak{a}=(a_1,\ldots,a_n)\subset A$, one lets
$R[T_1,\ldots,T_n]\rar \Rees_R(\mathfrak{a})=R[\mathfrak{a}\,T]$ be the graded map
sending $T_i$ to $a_i\,T$.
The {\em relation type} of $\mathfrak{a}$ is the largest degree of any
minimal system of homogeneous generators of the kernel $\cl J$.
Since the isomorphism $R[T_1,\ldots,T_m]/{\cl J}\simeq \Rees_R(\mathfrak{a})$ is graded, an application
of the Schanuel lemma to the graded pieces shows that the notion is independent of the set of generators
of $\mathfrak{a}$.

\begin{Corollary}\label{reltype}
Let $R$ be a Notherian ring and let $J\subset I$ be ideals in $R$ such that
\begin{enumerate}
\item[{\rm (i)}] $I/J$ has relation type at most $\ell$ as an ideal of $R/J$.
\item[{\rm (ii)}] $J\cap I^{t}=JI^{t-1}$ for every $t\leq \ell$.
\end{enumerate}
Then $\cl A _{R\surjects
  R/J}(I/J)\surjects \Rees_{R/J}(I/J)$ is an isomorphism.
\end{Corollary}
\demo
By Lemma~\ref{artinrees_isomorphism}, it suffices to show that $\ell$ is an upper bound for
the Artin--Rees number of $J$ relative to $I$.
Thus, Let $I$ be generated by elements $a_1,\ldots,a_n$
and let let $a\in J\cap I^t$, with $t\geq \ell$. Then there exists
 a homogeneous $F \in R[\TT]=R[T_1,\dots,T_n]$, of degree $t$, such that
$F(a_1,\ldots,a_n) = a$. Since $a\in J$, reducing modulo $J$ shows that $F$ is a
polynomial relation on $a_1,\ldots,a_n$ modulo $J$. Then, by assumption there are polynomials
$G_i,H_i\in R[\TT]$ of degrees $\ell,t-\ell$, respectively, such that $G_i(a_1,\ldots,a_n)\equiv 0\pmod{JR[\TT]}$
for $i=1,\ldots,n$ and $F\equiv\sum_iG_iH_i \pmod{JR[\TT]}$.
Write $ F= \sum G_iH_i +L$ for some homogeneous
polynomial $L\in JR[\TT]$ of degree $t$. Since
$$L(a_1,\ldots,a_n)\in JI^t\subset I^{t-\ell}(J\cap I^{\ell}),$$
then $G_i(a_1,\ldots,a_n)\subset J\cap I^{\ell}$
 and $H_i(a_1,\ldots,a_n)\in I^{t-\ell}$ for $t\geq \ell$. This shows that the element
$a=F(a_1,\ldots,a_n)\in I^{t-\ell}(J\cap I^{\ell})$, that is, $J\cap I^t=(J\cap I^{\ell})I^{t-\ell}$ for $t\geq \ell$.
\qed

\medskip

Using the notion of standard base we can add a tiny bit on the problem of describing
the associated primes of the Aluffi algebra.

\begin{Proposition}\label{torsion_primes}
Let $J\subset I\subsetneq R$ be ideals of a Noetherian ring $R$ such that
$I/J$ has a regular element.
If $p\in\ass_R(R/J)$ then
$\tilde{p}\in \ass_{\Rees_R(I)} (\cl A _{R\surjects R/J}(I/J))$.
\end{Proposition}
\demo
Let $p\in\ass _R(R/J)$.  Since $I/J$ contains a regular element,
one has  $I\not\subset p$. Let $a_1,\ldots,a_m$ be an
$I$-standard base of $p$, so that $p=(a_1,\ldots,a_m)$ and
$$\tilde{p}=\bigoplus_{t\geq 0}p\cap I^t=\bigoplus_{t\geq 0}\left(\sum_{i=1}^{m}a_iI^{t-\nu_i}\right),$$
where $\nu_i=\nu_I(a_i)$. Write $\nu=\max_i\nu_i$ for $i=1,\ldots,m$ and take $b\in I^{\nu-1}\setminus p$
(note that if $\nu=0$, which is a possibility, one means any $b\not\in p$).

Write $\mathfrak{R} =\bigoplus_{t\geq 0}J\cap I^t$ and $\mathfrak{A}=\bigoplus_{t\geq 0}JI^{t-1}$.

Say,  $p=J:a$. We claim that $\tilde{p}=\mathfrak{A}:ab$ which will prove that
$\tilde{p}$ is an associated prime of
$\cl A _{R\surjects R/J}(I/J))$ on $\Rees_R(I)$.

For this, let $c_t\in J\cap I^{t}\subset p\cap I^{t}=\sum_{i=1}^{m}a_iI^{t-\nu_i}$, with $t\geq 0$.
Then
$$b c_t \in \sum_{i=1}^{m}a_i\,bI^{t-\nu_i}\subset \sum_{i=1}^{m}a_iI^{t-\nu_i+\nu-1}
\subset \sum_{i=1}^{m}a_iI^{t-1}=
p\,I^{t-1}=(J:a)I^{t-1}\subset JI^{t-1}:a,$$
hence $b \tilde{p}\subset \mathfrak{A}:a $, hence $\tilde{p}\subset \mathfrak{A}:ab$.

For the inverse inclusion, since $\mathfrak{A}\subset \mathfrak{R} $, it follows that
$\mathfrak{A}:ab\subset \mathfrak{R} :ab=(\mathfrak{R} :a):b$.
Note that  $\widetilde{J:a}=(J:a)R[t]\cap\Rees_R(I)=(JR[t]:a)\cap\Rees_R(I)
=(JR[t]\cap\Rees_R(I)):a=\mathfrak{R} :a$.

Therefore $\mathfrak{A}:b a\subset \tilde{p}:b=\tilde{p}$.
Thus, $\tilde{p}=\mathfrak{A}:ab$ as was to be shown.
\qed

\subsection{Selected examples}

Let us agree to call a pair of ideals $J\subset I\subset R$ an {\em ${\cl A}$-torsionfree pair} if
the map $\cl A _{R\surjects A}(I/J)\surjects \Rees_{R/J}(I/J)$ is injective.

The examples we have in mind in this part are of the two sorts mentioned previously, namely,
of totally degenerate or non-degenerate pairs.
The first kind will be based on Proposition~\ref{residual_gens}.
For these, we let $R=k[\mathbf{X}]$ be an $\NN$-graded polynomial ring over a field $k$, $J\subset R$ is a homogeneous
ideal and $I\subset R$ is the Jacobian ideal of $J$, by which we always mean the ideal
$(J,\mathcal{I}_r(\Theta))$ where $r=\hht(J)$ and $\Theta$ stands for the Jacobian matrix of a set of generators
of $J$. One knows that this maybe a slack ideal, but it is well defined modulo $J$.

First we state a general format that implies an ${\cl A}$-torsionfree  pair $J\subset I$.

\begin{Example}\label{forms}\rm
Let $J\subset R=k[\xx]$ be an ideal  generated by forms of the same degree $d\geq 1$.
If $I=(J,\fm^r)$, where $\fm=(\xx)$ and $r\geq d$, then the pair $J\subset I$ is $\cl A$-torsionfree.
\end{Example}
To see this,  one uses Proposition~\ref{residual_gens}. Namely, it suffices to show that $J\cap\fm^{rt}\subset
J\fm^{r(t-1)}$ for every $t\geq 1$.
Let $a_1,\ldots,a_m$ be generators of $J$ of degree $d$ and let $F$ be a form in the $a_i$'s such that $F\in \fm^{rt}$.
Then $F=\sum _{i=1}^{m}G_ia_i$ where $G_i=\sum{a_{\underline{\alpha}}}\xx^{{\alpha}}\in R_{rt-d+\delta}$
for $\delta\geq 0$, since $R_{rt-d+\delta}=R_{r-d+\delta}\,R_{rt-r}$, so we can rewrite $G_i$ as
 $$G_i=\displaystyle\sum_{{|{\alpha}|=r-d+\delta}\atop {|{\beta}|=rt-r}}^{} \
 a_{{\alpha},{\beta}}\ \xx^{{\alpha}+{\beta}},\quad {\rm hence}\quad
F= \sum_{{|{\alpha}|=r-d+\delta}\atop\ }^{}\ \xx^{{\alpha}}
\left(\sum_{i=1\atop |{\beta}|=rt-r}^{s}(\xx^{{\beta}})a_i\right)$$

Therefore $F\in J\fm^{rt-r}= J\fm^{r(t-1)}$, as required.

\bigskip

Instances of this situation seem to be  any of the following ideals $J$ with respect to
the respective Jacobian ideal $I$.

\begin{enumerate}
  \item [{\rm (a)}] The defining ideal of the rational normal curve$\,${\rm ;}
  \item  [{\rm (b)}] The defining ideal of the Segre embedding of $\pp^r\times \pp^s$, with $r>1$ or $s>1\,${\rm ;}
  \item [{\rm (c)}] The defining ideal of the $2$--Veronese embedding of a projective space$\,${\rm ;}
   \item [{\rm (d)}] The ideal generated by the $4\times 4$ Pfaffians.
\end{enumerate}

It is well-known that  $J$ is the ideal of $2$--minors of the generic Hankel,
square, symmetric matrix,    respectively, and, lastly, the ideal generated by the maximal
Pfaffians of a $5\times 5$ skew-symmetric.

%Let $m$ denote the number of variables ${\bf x}$ in each case.
%In the first three cases, inspecting the nature of such matrices yields that the $2$--minors
% are typically of the form $x_i^2-x_lx_k$ and $x_ix_j-x_kx_l$,
% $0\leq i\neq j,k\neq l\leq m$.
% %This implies, in the three situations, that in the corresponding
%%Jacobian matrix $\Theta$ each row has at least $m-3$  zero entries.

Let $I$ denote the Jacobian ideal of $J$ on $R$.
Set $\fm=({\bf x})$ and write $\hht J=r\geq 2$.
We would need to prove that $I_r(\Theta)=\fm^r$, where $\Theta$ is the Jacobian matrix of
the natural generators of $J$.
A calculation of a good deal of cases gives evidence to this equality -- actually, it may be
the expression of a more general fact disguised under an inductive procedure.

Note that the case of the Segre embedding of $\pp^1\times \pp^1$ is exceptional, essentially because
it is a hypersurface. Here the ideal
of $1$-minors is $\xx$ which is of linear type, but clearly the defining ideal of the relative blowup
contains the equation of the dual to the (self-dual) quadric surface.

%\medskip

%The argument in the case of Pfaffians follows the same line of argument to show that $I_r(\Theta)=\fm^r$.
%Note that the ideal of the Pfaffians is the defining ideal of the Grassmanian ${\mathbb G}(1,4)$, hence its Jacobian
%ideal is a priori $(\xx)$-primary (generated in degree $r=3$).

\begin{Example}\rm Consider the monomial $x_1\cdots x_n\in R=k[x_1,\ldots,x_n]$
($n\geq 3$) and let $J\subset R$ be the ideal generated by its partial derivatives
$a_i=:x_1x_2\cdots\widehat{x_i}\cdots x_n$, for $i=1,...,n$.
If $I$ is the Jacobian ideal of $J$
the pair $J\subset I$ is $\cl A$-torsionfree.
\end{Example}
\demo Its well-known and easy that $J$ is a codimension $2$ perfect ideal
with Hilbert--Burch matrix

$$\phi=\left(
          \begin{array}{ccccc}
            x_1 & 0 & 0 &\ldots & 0 \\
            0 & 0 & 0 &  \ldots & x_2 \\
            \vdots & \vdots & \vdots  & \ldots & \vdots  \\
            0 & 0 & x_{n-2} & \ldots& 0 \\
            0 & x_{n-1} & 0 & \ldots & 0 \\
            -x_n & -x_n & -x_n &\ldots & -x_n \\
          \end{array}
        \right).
 $$

Setting
$\Delta_{i,j}:=\frac{\partial a_j}{\partial
x_i}=x_1x_2...\widehat{x_i}...\widehat{x_j}...x_n$,
and inspecting the Hessian matrix $\Theta$ of $x_1\cdots x_n$ -- a symmetric matrix --
one finds three basic types of $2\times 2$ minors, namely

\begin{itemize}
\item Principal minors:
$$ \det \left(
                                         \begin{array}{cc}
                                           0 &  \Delta_{i,j}\\
                                           \Delta_{i,j} & 0\\
                                         \end{array}
                                       \right)=\Delta^2_{i,j},
$$
one for each pair $i<j$;
\item vanishing minors:
$$\det \left(
         \begin{array}{cc}
           \Delta_{i,j} & \Delta_{i,j'} \\
           \Delta_{i',j} & \Delta_{i',j'} \\
         \end{array}
       \right)=0,$$
       one for each choice of row indices $1\leq i,i'\leq n$ and column indices $1\leq j,j'\leq n$;

\item       semidiagonal minors of typical form
$$ \Lambda_j:=\det\left(
                         \begin{array}{cc}
                           \Delta_{i,j} & 0 \\
                           * & \Delta_{i',i} \\
                         \end{array}
                       \right).$$
                       \end{itemize}
Since clearly, $\Lambda_j\in J$, we get that
the Jacobian ideal $I$ of $J$ is
 generated by $J$ and the squares of the second partial derivatives of $x_1\cdots x_n$, i.e.,
 $I=(J,\Delta^2_{i,j})$ for $1\leq i<j \leq n$.

As a side curiosity we note that
$I=(J,\mathcal{I}_2(\Theta))=(J,\mathcal{I}_{n-2}(\phi)^2)$, hence
$\sqrt{I}=\mathcal{I}_{n-2}(\phi)$, so in particular $I/J$ has codimension one.
This example will therefore yield a case of a height one ideal in $R/J$ which is $\cl A$-torsionfree,
but clearly not of linear type because its number of generators on $R/J$ is too large.

 Setting
$\Delta=(\Delta^2_{i,j}\ | \ 1\leq i<j \leq n)$, the usual algorithmic procedure
to find generators of the intersection of monomial ideals yields for any $t\geq 2$
$$J\cap \Delta^t=((x_i,x_j)\Delta_{i,j}^{2t}, (\mathfrak{F}))$$
 where $\mathfrak{F}$ is the set of all monomials in $\Delta^t$ excluding the monomials $\Delta_{i,j}^{2t}$ for $1\leq i<j \leq n$.
 Another calculation shows that  $(x_i,x_j)\Delta_{i,j}^{2t}\in J\Delta^{t-1}$, for $1\leq i<j \leq n$, and
 that $\mathfrak{F}\subset J^2\Delta^{t-2}$.
 This proves that $J\cap \Delta^t\subset JI^{t-1}$.
 \qed

 \begin{Question}\rm ($k$ algebraically closed) Let $J\subset R=k[x_1,\ldots, x_n]$ denote the homogeneous defining  ideal
of an arrangement of linear subspaces of dimension $n-3$ of $\mathbb{A}^n$.
If $I$ is the Jacobian ideal of $J$, when is $J\subset I$ an $\cl A$-torsionfree pair?
\end{Question}

Plausibly, a similar question can be posed about the Jacobian ideal of a hyperplane arrangement.

\begin{Example}\rm
Let $J\subset R=k[x_1,\ldots, x_n]$ denote the ideal of the coordinate points in
projective space ${\mathbb P}^{n-1}$, i.e.,
$J=(x_ix_j \ | \ 1\leq i<j\leq n)$.
If $I$ is the Jacobian ideal of $J$
the pair $J\subset I$ is $\cl A$-torsionfree.
\end{Example}
\demo Since $J$ contains all square-free monomials of degree $2$, it is rather transparent that
the Jacobian ideal $I$ of $J$ is generated by $J$ and pure powers of the variables.
Moreover, a closer inspection shows that, more precisely,
$$I=( J,\ x_1^{n-1},\ \ldots ,\ x_n^{n-1}).$$
Setting $\Delta:=(x_1^{n-1},\ \ldots ,\ x_n^{n-1})$, a procedure based on
finding generators of the intersection of monomial ideals yields for any $t\geq 2$
\medskip

 $$J \cap \Delta^t=(x_i^{t(n-1)}(x_1,\ldots,\widehat{x_i},\ldots, x_n),(\mathfrak{F})),$$
 where $(\mathfrak{F})$
is the set of all monomials in $\Delta^{t}$ excluding the monomials
$x_i^{t(n-1)}$ for $i=1,\cdots,n$. A calculation shows that
$x_i^{t(n-1)}(x_1,\ldots,\widehat{x_i},\ldots, x_n)\in J\Delta^{n-1}$  for $i=1,\cdots,n$ and that $(\mathfrak{F})\subset J^2\Delta^{t-2}$.
This proves that $J\cap \Delta^t\subset JI^{t-1}$.

\begin{Example}\label{four points} \rm
Let $J\subset R=k[x,y,z]$ denote the homogeneous defining ideal of the four points
$(1:0:0),(0:1:0),(0:0:1)$ and $(1:1:1)$ in the
projective plane $\mathbb{P}_k^2$ and let $I$ denote its Jacobian ideal. An easy calculation gives
$$J=(x^2-xz,y^2-yz)\quad\mbox{\rm and}\quad I=(x^2-xz,y^2-yz, x(2y-z),y(2x-z),(2x-z)(2y-z)).$$
In terms of the internal grading of the algebra, using the description in Proposition~\ref{torsion},
the torsion is generated by the appropriate residues of $\{xz^2(x-z),yz^2(y-z)\}\subset J\cap I^2$.
For further details see \cite[1.1]{Nejad}.
\end{Example}

The last two examples motivate the following

\begin{Question}\rm ($k$ algebraically closed) Let $J\subset R=k[x_1,\ldots, x_n]$ denote a radical homogeneous ideal
of codimension $n-1$ (i.e., the ideal of a reduced set of points).
If $I$ is the Jacobian ideal of $J$, when is $J\subset I$ an $\cl A$-torsionfree pair?
\end{Question}

\section{The gradient  Aluffi  algebra of a projective hypersurface}\label{gradients}

In this section we will deal with the case where $J$ is generated by the equation of a reduced hypersurface.

Thus, let $J=(f)\subset R$ be a principal ideal,
where $R=k[X_1,\ldots,X_n]$.
We will focus on the Jacobian ideal $I_f=(f, \partial f/\partial X_1, \ldots, \partial f/\partial X_n)$.
We are particularly motivated by the problem as to when $I_f$ is an ideal of linear type.
In general, if $f$ is not Eulerian, the local number of generators of $I_f$ maybe an early obstruction
to this property.
We will consider the case where $f$ is homogeneous in the standard grading of the polynomial ring and
its degree is not a multiple of the characteristic -- hence, $f\in I_f$. In this context the ideal $I_f$ will often be of linear type.

We call $I_f$ the {\em gradient} ideal of $f$ and the corresponding algebra $ \cl A _{R\surjects R/(f)}(I_f/(f))$,
the {\em gradient Aluffi algebra} of $f$.
By Proposition~\ref{dimension_bounds} (d), it  is equidimensional
of dimension $\dim R=n$.

\subsection{Preliminaries on the gradient ideal}

If $f$ defines a smooth hypersurface then $I_f$ is irrelevant, i.e., is generated by a regular sequence, hence is of linear type.
We regard this case as uninteresting and assume that $f$ has singularities. This entails $\hht (I_f/(f))\leq n-2$.
If moreover $f$ is reduced then $\hht (I_f/(f))\geq 1$.
For $n=3$ we therefore find $\hht (I_f/(f))=1$.
Ideals of height $1$ in non-regular rings of dimension $2$ are  a tall order and
typically involve a non-trivial primary decomposition.

Thus, even over projective plane curves the structure of the gradient Aluffi  algebra seems to be fairly intricate.
Note that for $n=3$, the ideal $I_f$ is an almost complete intersection.
Since we regard the linear type case as a limit situation we would like to understand this case first.

Now, for an  almost complete intersection  this property is fairly under control.
For convenience we file the following general result, which collects in a more detailed version several known
facts about an almost complete intersection (see \cite[Proposition 3.7]{conormal}, also
\cite[Proposition 8.4, Proposition 10.4, Remark 10.5]{Trento}).

\begin{Lemma}\label{aci}
Let $R$ denote a Cohen--Macaulay local ring and
let $I\subset R$ denote a proper ideal of height $h\geq 0$. Assume that
\begin{itemize}
\item $I$ is a strict almost complete intersection
{\rm (}i.e., minimally generated by $h+1$ elements{\rm )}
\item $R/I$ is equidimensional  {\rm (}i.e., $\dim R/I=\dim R/P$
for every minimal prime $P$ of $R/I${\rm )}
\item $I$ satisfies the so-called sliding depth inequality $\depth R/I\geq \dim R/I-1$.
\end{itemize}

Let $R^m\stackrel{\phi}{\lar} R^{h+1} \lar I \lar 0$ stand for a minimal free presentation of $I$ as an $R$-module.
The following conditions are equivalent:
\begin{enumerate}
\item[{\rm (1)}] $\hht I_1(\phi)\geq \hht I +1$
\item[{\rm (2)}] $I_P$ is a complete intersection for every minimal prime $P$ of $R/I$
\item[{\rm (3)}] $I$ is of linear type.
\end{enumerate}
\end{Lemma}
\demo
(1) $\Rightarrow$ (2) Localizing at such a prime leaves some element of $I_1(\phi)$ invertible,
so up to an elementary transformation on $\phi_P$ the local presentation has the form
$$R_P^{m-1}\oplus R_P\stackrel{\phi_P}{\lar} R_P^h\oplus R_P \lar I_P \lar 0,$$
with
\arraycolsep=10pt

$$\phi_P=\left(
\begin{array}{c@{\quad\vrule\quad}c}
\raise5pt\hbox{$1$}&\raise5pt\hbox{$0$}\\[-6pt]
\multispan2\hrulefill\\%[-6pt]
\lower2pt\hbox{$0$}&\lower2pt\hbox{$\psi$}\\%[-2pt]
\end{array}
\right)
$$
Therefore, we get a free presentation $R_P^{m-1}\stackrel{\psi}{\lar} R_P^h \lar I_P \lar 0$, thus
showing that $I_P$ is generated by (at most) $h$ elements.

(2) $\Rightarrow$ (3) By \cite[Proposition 10.4]{Trento} the symmetric algebra of $I$ is a Cohen--Macaulay ring.
 Therefore, by \cite[Proposition 8.4]{Trento} it suffices to show that $\hht I_Q\leq \hht Q$ for every prime $Q\subset R$
 containing $I$. Let $P$ be a minimal prime of $R/I$ contained in $Q$. If $P=Q$ the hypothesis
 guarantees the inequality. Otherwise
$\hht(Q)\geq h+1$. But $\hht(I_Q)= \hht I_P=h$ because $R/I$ is equidimensional, hence we are through.

(3) $\Rightarrow$ (1) By \cite[Lemma 8.2 and Proposition 8.4]{Trento} one has $\hht I_t(\phi)\geq \rk (\phi)-t+2$
 for every $1\leq t\leq \rk (\phi)=h$. In particular, $\hht I_1(\phi)\geq h-1+2=\hht I +1$.

\qed

\begin{Corollary}\label{aci_dim3}
Let $f\in R=k[x_1,\ldots,x_n]$ stand for a reduced homogeneous polynomial.
Assume that the singular locus of $V(f)\subset \pp^{n-1}$ consists of a nonempty set of
points.
The following are equivalent:
\begin{enumerate}
\item[{\rm (1)}] The coordinates of the vector fields of
${\mathbb P}^{n-1}$ vanishing on $f$ generate an irrelevant ideal.
\item[{\rm (2)}] Locally at each singular point
of $V(f)$ the gradient ideal is a complete intersection.
\item[{\rm (3)}] The gradient ideal of $f$ is of linear type.
\end{enumerate}
\end{Corollary}
\demo A vector field $v=\sum_{i=1}^n a_i\partial /\partial x_i$ vanishes on $f$
if and only if  $\sum_{i=1}^n a_i\partial f/\partial x_i=0$.
Therefore the coordinates of all such vector fields generate the ideal of $1$-minors of a syzygy matrix
of the gradient ideal. The result then follows  from Lemma~\ref{aci} once its hypotheses
are verified in this setup, as we next proceed to see.

Since $f$ is  assumed to be reduced, whose singular locus is a nonempty set of isolated singularities,
its gradient ideal is a (homogeneous) ideal of codimension $n-1$, hence can only have minimal primes
of codimension $n-1$ in $R$. Therefore it is equidimensional.

Finally, the depth condition is trivially verified for the numbers in question.
\qed

\bigskip

So much for the linear type property.
Clearly, this property implies that the partial derivatives are algebraically independent over $k$.
The latter property in turn reads geometrically to the effect that the polar map associated to the
hypersurface $V(f)\subset \pp^{n-1}$ is dominant. In characteristic zero this is tantamount to saying
that the Hessian of $f$ does not vanish (cf. \cite{CRS} for a detailed account on this).

\medskip

The following result collects parts of the main backstage for the Aluffi gradient algebra.

\begin{Theorem}\label{central_aluffi_curve}
Let $k$ denote an infinite field, let $f\in R=k[\xx]=k[x_1,\ldots,x_n]$ be a reduced homogeneous polynomial
whose degree is not a multiple
of the characteristic of $k$ and let $I_f\subset R$ denote the corresponding gradient ideal.
Assume that
\begin{enumerate}
\item[{\rm (i)}] The singular locus of $V(f)\subset \pp^{n-1}$ consists of a nonempty set of
points {\rm (}equivalently, $\dim R/I_f=1${\rm )}
\item[{\rm (ii)}] The partial derivatives of $f$ are algebraically
independent over $k$.
\end{enumerate}
Then:
\begin{enumerate}
\item[{\rm (a)}] The minimal primes of the gradient Aluffi algebra on $\Rees_R(I_f)$ are
\begin{itemize}
\item The minimal prime ideals of $\Rees_{R/(f)}(I_f/(f))$, all of the form $\sum_{t\geq 0}(p)\cap I^t$
for a prime factor $p$ of $f$
\item The extended ideal $(\xx)\Rees_R(I_f)$
\item Prime ideals whose lifting to $R[\TT]=R[T_1,\ldots,T_n]$ from a presentation $R[\TT]/{\cl A}\simeq \Rees_R(I)$
have the form $(P, {\mathfrak f})$, where $P\subset R$ is a minimal prime of $R/I_f$ and ${\mathfrak f}$
is an irreducible homogeneous polynomial in $k[\TT]$.
\end{itemize}
\item[{\rm (b)}] The gradient Aluffi  algebra has non-trivial torsion.
\item[{\rm (c)}] $I_f$ is an ideal of linear type {\rm (}respectively, weakly of linear type{\rm )}
if and only if the natural surjection
$${\cal S}_{R/(f)}({I_f/(f)})\surjects \cl A _{R\surjects R/(f)}(I_f/(f))$$
is an isomorphism {\rm (}respectively, an isomorphism in all high degrees{\rm )}.
\item[{\rm (d)}] The symmetric algebra ${\cal S}_{R/(f)}({I_f/(f)})$ is Cohen--Macaulay$\,${\rm ;} in particular, if
$I_f$ is of linear type then the gradient Aluffi algebra is Cohen--Macaulay.
\end{enumerate}
\end{Theorem}
\demo
(a) We apply Proposition~\ref{minimal_primes}, from which the first set of minimal primes is
clear.

To see that $(\xx)\Rees_R(I_f)$ is a minimal prime thereof one proceeds as follows.
There is a presentation of the gradient Aluffi  algebra
\begin{equation}\label{presenting_aluffi_gradient}
\cl A _{R\surjects R/(f)}(I_f/(f))\simeq R[\TT]/({\cl J}_f, f, \sum_{i=1}^n x_iT_i),
\end{equation}
where ${\cl J}_f$ denotes the defining ideal of the Rees algebra $\Rees _R(I_f)$ on $R[\TT]$.

Since the partials are homogeneous of the same degree  algebraic independence over $k$
is tantamount to analytic independence (i.e., the relations of the
generators of $I_f$ have coefficients in the ideal $(\xx)$).
Therefore, the result follows from
Proposition~\ref{lying_deep}, (ii).

Let ${\cl P}$ be a minimal prime of $\cl A _{R\surjects R/(f)}(I_f/(f))$ whose
contraction $P$  to $R$ contains $I_f$ and is properly contained in $(\xx)$.
By Proposition~\ref{minimal_primes}, ${\cl P}=(P,{\cl P}_+)$.
By assumption (i) it follows that $P$ is a minimal prime of $R/I_f$, hence has height $n-1$.
By Theorem~\ref{hypersurface}, $\cl A _{R\surjects R/(f)}(I_f/(f))$ is equidimensional.
Therefore the lifting of ${\cl P}$ to $R[\TT]$ has height $n$.
Since the lifting of any minimal generator of $({\cl P}_+)$ is irreducible in $k[\TT]$
it follows  that the lifting of ${\cl P}$ to $R[\TT]$ has the required form.

\medskip

(b)
The defining equation of the dual curve to $f$ belongs to the presentation ideal of $\Rees_{R/(f)}(I_f/(f))$
on $R[\TT]$  and not to $(\xx)R[\TT]$, hence by (a) it does not belong to the defining ideal
of $\cl A _{R\surjects R/(f)}(I_f/(f))$ on $R[\TT]$.

\medskip

(c)
This is a straightforward application of Proposition~\ref{lineartype}.
The argument for the weak version of the property of being of linear type is exactly the same
as in [loc.cit.].

\medskip

(d) We apply the criterion of \cite[Theorem 10.1]{Trento}. Namely, we have to verify the
following conditions:
\begin{enumerate}
  \item [{\rm (A)}] $\mu(I_f/(f)_{P/f})\leq \hht (P/(f))+1=\hht P$, for every prime ideal $P\supset I_f$ of $R$.
   \item [{\rm (B)}] $\depth (H_i)_{P/(f)}\geq \hht (P/(f))-\mu(I_f/(f)_{P/(f)})+i=\hht P-\mu(I_f/(f)_{P/(f)})+i-1$,
   for every prime ideal $P\supset I_f$ of $R$ and every $i$ such that $0\leq i\leq  \mu(I_f/(f)_{P/(f)})-\hht (I_f/(f)_{P/(f)})$,
   and where $H_i$ denotes the $i$th Koszul homology module of the partial derivatives on $R/(f)$.
\end{enumerate}
Note that the primes containing $I_f$ are
$\fm=(x_1,\ldots,x_n)$ and its minimal primes, the latter all  of height $n-1$.

(A) Since $I_f$ itself is generated by $n$ elements, it suffices to check the minimal primes.
Thus, let $P\subset R$ be such a prime. Say, without lost of generality,
that $x_n\not\in P$. Because of the Euler relation,  $\partial{f}/\partial{x_n}\in I_f$
locally at $P$ and module $(f)$. Therefore,
locally at $P$ and module $(f)$, $I_f$ is generated by $n-1=\hht (P)$ elements.

(B) If $P$ is a minimal prime of $I_f$ we saw in (A) that $\mu(I_f/(f)_{P/(f)})=n-1$.
Since $\hht P=n-1$, the condition is trivially verified as $i=0,1$.

Thus, let $P=\fm$.
Again, an easy inspection of the numbers tell us that only the case where $i=2$ needs an argument and, in this case,
one has to prove that $\depth (H_2)_{\fm/(f)}\geq 1$.
Localize $R$ at $\fm$ and update the notation, so $R:=R_{\fm}\supset I_f:={I_f}_{\fm}\supset (f)=(f)_{\fm}$
and $H_2:=(H_2)_{P/(f)}$.

Now, $f$ is a nonzero element in $I_f$ and $I_f$ has grade $n-1$ in $R$.
Thus, there is a regular sequence in $I_f$ of length $n-1$ starting with $f$.
Write $L\subset I_f$ for the ideal generated by this regular sequence.

The following isomorphism is well know (see, e.g., \cite[Theorem 1.6.16]{BH}):
$$H_2\simeq \hom_{R/(f)}\left(\frac{R/(f)}{I_f/(f)},\frac{R/(f)}{L/(f)}\right)\simeq
\hom_{R/(f)}\left({R}/{I_f},{R}/{L}\right).$$

Therefore $\ass_{R/(f)}(H_2)={\rm Supp}_{R/(f)}( R/I_f)\cap\ass_{R/(f)}(R/L)\subset \ass_{R/(f)}(R/L)$.
But $L$ is generated by a regular sequence of length $n-2$ modulo $f$ by construction, while $\dim R/(f)=n-1$.
It follows that $\fm/(f)\not\in \ass_{R/(f)}(R/(J,f))$, hence
$\fm/(f)\not\in \ass_{R/(f)}(H_2)$.
\qed

\begin{Remark}\rm
 For $n\leq 4$, if  the partial derivatives are $k$-linearly independent  then the result of Gordan--Hesse--Noether  implies that they are
algebraically independent over $k$ (see \cite[Proposition 2.7]{CRS} for a proof of the case $n=3$
based on an observation of Zak).
Thus, the assumption in this range is just linear independence.

As to (c), it's valid
with no restriction when $f$ is reduced since the defining equations of the dual variety to the hypersurface $V(f)$
belong to the presentation ideal ${\cal A}_f$ and, moreover, the ideal generated by these contains properly the defining ideal
 of the polar map of $V(f)$ (see \cite[Remark 2.4]{CRS}).
\end{Remark}

\begin{Example}\rm
Here is a simple illustration.
Let $f=x^2y^2+x^2z^2+y^2z^2$, the equation of a plane quartic with $3$ ordinary nodes.
The minimal primes of the corresponding gradient Aluffi algebra, lifted to $k[x,y,z,T,U,V]$, are $(x,y,z)k[x,y,z,T,U,V]$, $(x,y,V)$, $(x,z,U)$, $(y,z,T)$ and its lifted torsion.
Since $I_f$ is of linear type (see next section), these are of course the minimal primes of the
symmetric algebra ${\cl S}_{R/(f)}(I_f/(f))$.
\end{Example}

\subsection{Gradient ideals of linear type along a family}

In general, the gradient ideal $I_f$ will not be of linear type. This subsection will prepare the ground
to showing that if $f$ is the equation of an irreducible plane rational {\em quartic} then $I_f$ is an ideal of linear type.

We can immediately show simple cases of rational plane quintics whose corresponding gradient ideals are not of linear type.
Moreover, though the corresponding gradient Aluffi algebras are equidimensional, they tend to behave quite erratically from the viewpoint
of the Cohen--Macaulay locus and of the associated primes.
It is apparent that this behavior reflects the nature of the singularities, but it is in general quite misterious.

\begin{Example}\label{bad_quintics}\rm
Let $f=y^4z+x^5+x^3y^2$. Then $I_f=(x^2(5x^2+3y^2), y(2x^3+4y^2z), y^4)$.
Canceling the common factor among the last two generators, gives rise to the obvious Koszul relation.
From this it immediately follows that the radical of the ideal generated by the coordinates of the syzygies of
$I_f$ has $x,y$ among its minimal generators.
The rest follows by inspection, as it is not difficult to verify that no syzygy coordinate has as term a pure $z$-power.
By Corollary~\ref{aci_dim3}, $I_f$ is not of linear type.

Of course everything in this example is easily obtained by machine computation.
The three algebras $\cl S _{R/(f)}(I_f/(f))\surjects \cl A _{R\surjects R/(f)}(I_f/(f))   \surjects \Rees _{R/(f)}(I_f/(f))$
are all distinct, but of the same dimension. The leftmost is Cohen--Macaulay, while the Aluffi algebra has no embedded primes
though it is not Cohen--Macaulay.

Now let $f=zy^2(x^2+y^2)+x^5+y^5+x^3y^2$.
Here the symmetric algebra is Cohen--Macaulay, while
 the Aluffi algebra has embedded primes.
\end{Example}

In this part we study families of singular plane curves and a certain ``partial'' gradient ideal for
the linear type property and the corresponding Aluffi algebra.
We start by making clear what we mean by a family for our purposes.
Note that the considerations that follow work {\em ipsis litteris} for hypersurfaces.

Let $k[\uu]=k[u_1,\ldots,u_m]$ stand for a polynomial ring over the field $k$ and let $F\in S:=k[\uu][x,y,z]$
denote a polynomial which is a form on $x,y,z$. We give $S$ the structure of standard graded ring over $k[\uu]$.
The basic assumption is that the content of $F$ with respect to the $\uu$-coefficients
is $1$. Then $F$ is a non-zero-divisor on $k[\uu]/I$ for every ideal $I\subset k[\uu]$, hence
$\mbox{\rm Tor}_{k[\uu]}^1(k[\uu]/I, k[\uu][x,y,z]/(F))=\{0\}$ for any such ideal.
This gives that the inclusion $k[\uu]\subset k[\uu][x,y,z]/(F)$ is flat, hence defines a family of curves in $\pp^2$ over
the parameters $\uu$.

Thus, we speak of a {\em family of plane curves} over the {\em parameters} $\uu$ when referring to this setup.
We will of course adhere to the terminology of calling {\em general member} of the family the equation of the plane curve
obtained by substituting general values in $k$ for $\uu$.
Moreover, our interest lies on the case where the general member of the family is a reduced singular plane curve.
In this case we speak of a {\em family of plane singular curves}.

\medskip

In the sequel we will assume moreover that $m\leq {{d+2}\choose {2}}-1$, where $d$ is the (homogeneous) degree of $F$ in $x,y,z$
and that $F$ has the form
\begin{equation}\label{F_form}
F=\phi_0(x,y,z)+\sum_{j=1}^m u_j\phi_j(x,y,z),
\end{equation}
where $\{\phi_j(x,y,z)\,|\, 0\leq j\leq m\}$ is a set of monomials of degree $d$ in $x,y,z$, and $\phi_0(x,y,z)\neq 0$.

\bigskip

Note that the form of $F$ depends on the singular points of the general member.
Thus, it makes sense to speak about a {\em normal form} or {\em canonical form} of $F$ depending on this singular locus.
Our convention is that such a normal form is to be obtained through projective transformations applied to the $x,y,z$-coordinates
allowing coefficients from $k[\uu]$. Besides, in order to account for
{\em degeneration of singularities} of the general member we need correspondingly to consider certain {\em degeneration ideals} in
the parameter ring $k[\uu]$.

Write
$$F\equiv \phi_0(x,y,z)+\psi(x,y,z,u_1,\ldots,u_m),$$
as in (\ref{F_form}),  where  $\phi_0(x,y,z)$ involves the singularity type in terms of the projectivized  tangent cones on suitable affine pieces.

\begin{Example}\rm Let us write a normal form for the  family
of irreducible singular quartic plane curves such that the singular
locus of the general member consists of one simple node - note that at
this point it is not totally clear that there exists at all such a family
 in the sense we  established, since we must first obtain some $F\in S$ that works.
  By projectivities one can assume that the node is $P=(0:0:1)$ and the tangent
  cone at $z\neq 0$ has equation $xy$. Since the general
member ought to vanish at $P$ then we may omit the terms in $z^4, z^3x$ and $z^3y$.
Thus, an intermediate step towards a normal form is
$$F=xyz^2+u_1x^3z+u_2x^2yz+u_3xy^2z+u_4y^3z+u_5x^4+u_6x^3y+u_7x^2y^2+u_8xy^3+u_9y^4.$$
We can see that the specialization of $F$ by $k$-values factors properly if
both $u_1$ and $u_5$ have vanishing $k$-values; similarly,
if both $u_4$ and $u_9$ have vanishing $k$-values.
Thus, for writing a normal form we may incorporate $x^4$ and $y^4$
as terms of $\phi_0(x,y,z)$.
Finally, the projectivity $x=x,y=y,z=z-\frac{1}{2}(u_2x+u_3y)$
(characteristic $\neq 2$) allows to eliminate the terms in $x^2yz$ and $xy^2z$.
Up to renaming parameters, this yields the following normal form:
$$F=xyz^2+x^4+y^4+u_1x^3z+u_2y^3z+u_3x^3y+u_4x^2y^2+u_5xy^3$$
\end{Example}

\subsection{Degeneration for the linear type condition}

A piece of difficulty regarding the question as to how the property of being of linear type moves on a family is that
this property is neither kept by
specialization nor by generization.
This difficulty permeates  the theory by often conflicting with the usual degeneration conditions
considered in the realm of families of hypersurfaces.

The normal form has degenerations to other normal forms whose general member has more involved singularities or even acquires new
singular points.
The following example may illuminate this phenomenon.

\begin{Example}\rm Consider the family of irreducible rational plane
quartics with exactly three nodes. In \cite[Lemma 11.3]{Gibson} a normal
form is given of a family whose general member is an irreducible quartic with three double points, namely
$$F=\lambda x^2y^2+\mu x^2z^2+\nu y^2z^2+2xyz(u_1 x+u_2 y+u_3 z),\; \lambda\nu\mu\neq 0.$$
To get a normal form whose general member is an irreducible quartic with three nodes, substituting $x=(\nu/\lambda\mu)^{1/4}x,\ y=(\mu/\lambda\nu)^{1/4}y,\ z=(\lambda/\mu\nu)^{1/4}z $ and renaming,
one obtains the normal form
$$F=x^2y^2+x^2z^2+y^2z^2+2xyz(u_1 x+u_2 y+u_3 z).$$
\end{Example}

Note that for $k$-values $u_1=\pm 1$, one of the nodes degenerates into a cusp and, similarly, for $u_2=\pm 1$ or $u_3=\pm 1$.
Thus, the general member requires that the $k$-values of the triple $(u_1,u_2,u_3)$ do not lie on the hypersurface $V((u_1^2-1)(u_2^2-1)(u_3^2-1))$
in order that it have exactly three nodes.

Requiring that the general member acquire no new singular points besides the three nodes imposes yet another obstruction.
Of course, in the present low degree $4$, because of genus reason there will be new singular points only if the general member properly factors.
As we will see this obstruction is precisely given by the hypersurface  whose equation is the discriminant $2u_1u_2u_3+u_1^2+u_2^2+u_3^2-1$ of a
suitable conic (see Section~\ref{rational_quartics_all}).

The following is a basic result for this part.
It would mostly suffice for it to assume that the $\uu$-coefficients of the terms of $F$ be algebraically independent
over $k$.
We observe that a similar result holds for families of hypersurfaces whose general member is  reduced and irreducible and,
moreover, the singular locus is a nonempty set of points.

\begin{Theorem}\label{degeneration_vs_lineartype}
Let $F$ denote a family of singular plane curves of degree $d\geq 2$, on parameters ${\mathbf a}=u_1,...,u_m$, whose general member is reduced
and irreducible.
Write $S=k[{\uu}][x,y,z]$. Let $I_F\subset S$ denote the ideal generated by the $x,y,z$-partial derivatives of $F$ and
let $\mathcal{I}\in S$ stand for the ideal of $1$--minors of the syzygy matrix of $I_F$. Then:
\begin{itemize}
  \item [{\rm (a)}] $I_F$ has codimension $2$
  \item [{\rm (b)}] $\mathcal{I}$ has codimension at most $3$
  \item [{\rm (c)}] If $k$ is algebraically closed of characteristic zero, the following are equivalent:
  \begin{itemize}
    \item [{\rm (i)}] $\mathcal{I}$  has codimension $3$.
    \item [{\rm (ii)}]  The contraction of the ideal  $\mathcal{I}:(x,y,z)S\,^{\infty}$ to $k[{\uu}]$ has codimension $\geq 1$.
     \item [{\rm (iii)}] The plane projective curve $F(\alpha)\in k[x,y,z]$ obtained by evaluating $\uu\mapsto \alpha$
     off a set of  codimension $\geq 1$ in  $\mathbb{A}^m_k$
      has gradient ideal of linear type.
      \item [{\rm (iv)}]  There is some $\alpha\in \mathbb{A}^m_k$ for which the evaluated ideal $\mathcal{I}(\alpha)\in k[x,y,z]$
      has codimension $3$.
  \end{itemize}
\end{itemize}
\end{Theorem}
\demo
(a) Clearly, $\codim (I_F)\leq 3$. Since the general member is singular and reduced its gradient ideal
has codimension $2$.
This forces $\codim (I_F)=2$.

\medskip

(b) We go more algebraic: the ring $S=k[{\uu}][x,y,z]$ is standard graded with
$S_0= k[{\uu}]$.
Since $F$ is a homogeneous polynomial in this grading,  its partial derivatives are homogeneous of same degree.
We claim that any syzygy of the partial derivatives has coefficients in $(x,y,z)S$.
Indeed, since the partial derivatives are homogeneous of the same degree $d\geq 1$ in $x,y,z$,
any syzygy is homogeneous of non-negative degree in $x,y,z$.
Now, if  a syzygy would happen to be of degree $0$,
i.e., with all its coordinates in the zero degree part $k[{\uu}]$,
this would force, by reading the relation in degree $0$,  a polynomial relation among the coefficients of degree $0$ of the partials, hence a
polynomial relation of ${\uu}$, which is nonsense since these are indeterminates over $k$.

Incidentally, note this argument breaks down for  syzygies of a higher order as the first syzygies may have different
degrees in $S$.

\medskip

(c)
(i) $\Rightarrow$ (ii)
Write $\fm=(x,y,z)S$. Since  $\codim (\mathcal{I}) = 3$ and $\fm$ is
a minimal prime therein by the proof of (b),  the saturation $\mathcal{I}:\fm^{\infty}$ picks up the primary
components of $\mathcal{I}$ not containing $\fm$. This shows that $(\mathcal{I}:\fm\,^{\infty})\cap k[{\uu}]\neq \{0\}$.

\medskip

(ii) $\Rightarrow$ (iii)  Let $g=g(\uu)\in ({\mathcal{I}}:\fm^{\infty})\cap k[\uu]$ be any nonzero element.
By hypothesis, $g$ conducts a power of $(x,y,z)S$ inside $\mathcal{I}$.
Giving $\uu$ $k$-values $\alpha$ off $V(({\mathcal{I}}:\fm^{\infty})\cap k[\uu])$ yields a power of the maximal ideal
$(x,y,z)\subset k[x,y,z]$
inside the image $\mathcal{I}(\alpha)$ of $\mathcal{I}$ by this evaluation. Let $f=F(\alpha)\in k[x,y,z]$ denote the
member of the family thus obtained. Then $\mathcal{I}(\alpha)\subset I_1(\phi)$, where $\phi$ denotes
the syzygy matrix of the partial derivatives of $f$.
This shows that $I_1(\phi)$ is $(x,y,z)$-primary.
Therefore, the result follows from Corollary~\ref{aci_dim3}.

\medskip

(iii) $\Rightarrow$ (iv) The hypothesis is that the gradient ideal of the general member of the family is of linear type.
Again by Corollary~\ref{aci_dim3} this implies that the ideal of $1$-minors of such a plane curve
has codimension $3$. On the other hand, for general value $\alpha$ of $\uu$, one has
$\mathcal{I}(\alpha)=I_1(\phi_{\alpha})$, where $\phi_{\alpha}$ stands for the syzygy matrix of $F(\alpha)$.

\medskip

(iv) $\Rightarrow$ (i)
 By definition,  $\mathcal{I}(\alpha)=(\mathcal{I},\uu-\alpha)/(\uu-\alpha)$ upon identifying $k[\uu,x,y,z]/(\uu-\alpha)=k[x,y,z]$
under the surjection $k[\uu][x,y,z]\surjects k[x,y,z]$ such that $\uu\mapsto \alpha$.
Since $\uu-\alpha$ is a regular sequence in a polynomial ring we easily get
\begin{eqnarray*}%\label{special_to_general}
&&\hht ((\mathcal{I},\uu-\alpha)/(\uu-\alpha))=\hht (\mathcal{I},\uu-\alpha)-\hht(\uu-\alpha)\leq\\
&& \hht(\mathcal{I})+\hht(\uu-\alpha)-\hht(\uu-\alpha)=\hht(\mathcal{I}),
\end{eqnarray*}
which implies the result.
\qed

The following example illustrates the various obstructions.

\begin{Example}\rm
The one-parameter family $F=y^4z+x^5+ux^3y^2$ (see Example~\ref{bad_quintics}) is such that $\mathcal{I}$ has codimension $2$,
hence the gradient ideal of the  general member $F(\alpha)$ of the family is not of linear type.
Clearly, it follows that $\mathcal{I}({\alpha})$ has height $\leq 2$ for any $\alpha\in k$.
A computation with {\em Macaulay} gives moreover that the associated primes of $S/\mathcal{I}$ are $(x,y)\subset (x,y,z)$.
Perhaps surprisingly, the
special member $F({\bf 0})$ is easily seen to have gradient ideal of linear type, i.e., the ideal of $1$-minors of the syzygy matrix
of the special member $F({\bf 0})$, obtained by evaluating $F$ at ${\bf 0}$, has codimension $3$.
This simple example shows that
the property in question does not deform to the generic member.
\end{Example}

Under the equivalent assumptions of item (c) of Theorem\ref{degeneration_vs_lineartype}, one can give
the approximate structure of the contracted ideal in item (ii).

\begin{Proposition}\label{when_cod_is_one}
Let the assumptions be those of Theorem~\ref{degeneration_vs_lineartype} and
assume that $\mathcal{I}$  has codimension $3$.
If $\mathcal{I}$ has a minimal prime of codimension $3$ other than $\fm=(x,y,z)S$
then $(\mathcal{I}:\fm\,^{\infty})\cap k[{\uu}]$ has codimension $1$.
If, moreover, $\mathcal{I}$ is pure-dimensional then this contraction is
a principal ideal.
\end{Proposition}
\demo Let $\mathfrak{p}$ be a minimal prime of codimension $3$ of $\mathcal{I}$
other than $\fm$.
Since $I_F\subset \mathcal{I}$ (because of the Koszul relations) then $\mathfrak{p}$
contains a minimal prime of $I_F$. But the latter are of two sorts: either the extensions of the
minimal primes (in $k[x,y,z]$) of the singular points of the general curve
of the family, or else minimal primes of codimension $3$. In the first case, $\mathfrak{p}$  contains two independent $1$-forms
in $k[x,y,z]$ which, up to a projective change of coordinates, can be assumed to be $x,y$.
Clearly, these forms are part of a minimal set of generators of $\mathfrak{p}$
and $(x,y,\mathcal{I})\subset \mathfrak{p}$.
But $\mathcal{I}$ is homogeneous in the variables $k[x,y,z]$, generated in
positive such degrees (see the argument in the proof of (b)), hence there exist suitable
polynomials $g_j(\uu)\in k[\uu]$ and integers $k_j\geq 1$, $1\leq j\leq s$, such that
$$(x,y,\mathcal{I})=(x,y,\,z^{k_j}g_j(\uu),1\leq j\leq s).$$
Since $z\notin \mathfrak{p}$, necessarily some $g_j(\uu)\in \mathfrak{p}$.
But then $\mathfrak{p}=(x,y,p(\uu))$, for some prime factor of $g_j(\uu)$.

We now deal with the case where $\mathfrak{p}$ contains  a minimal
prime of $I_F$ of codimension $3$, hence coincides with it.
Let $\mathfrak{p}_i, 1\leq i\leq r$ denote the minimal primes of codimension $3$
of $\mathcal{I}$. Each of these, by the previous argument, has a minimal generator
$p_i(\uu)\in k[\uu]$.
Then
$$\left(\prod_i\mathfrak{p}_i\right)\cap k[\uu]=(\prod p_i(\uu)),$$
a principal ideal.
On the other hand,
$$\sqrt{\mathcal{I}:\fm^{\infty}}\subset \bigcap \sqrt{\mathcal{P}_i:\fm^{\infty}}
\subset \bigcap (\sqrt{\mathcal{P}_i}:\fm^{\infty})=\bigcap \sqrt{\mathcal{P}_i}=
\bigcap \mathfrak{p}_i=\sqrt{\prod_i\mathfrak{p}_i},
$$
where $\mathcal{P}_i$ denotes the $\mathfrak{p}_i$th primary component of $\mathcal{I}$.
This proves that the contraction $(\mathcal{I}:\fm\,^{\infty})\cap k[{\uu}]$
has codimension $\leq 1$, hence is exactly $1$ by Theorem~\ref{degeneration_vs_lineartype}.

The additional assertion at the end of the statement is now clear.
\qed

\subsection{Rational quartics}\label{rational_quartics_all}

We review some preliminaries about rational quartics, the basic reference being
\cite{Wall}.

An irreducible rational quartic having only
 double points can be obtained as a rational transform from a
non-degenerate conic by means of one of the three basic plane quadratic Cremona maps:

\begin{enumerate}
\item[{\rm (1)}] $\pp^2\dasharrow \pp^2$ with defining coordinates $(yz:zx:xy)$

The base locus of this Cremona map consists of the points $(1:0:0)$, $(0:1:0)$ and $(0:0:1)$,
each with multiplicity one (in the classical terminology, three proper points -- see \cite{Alberich}).

\item[{\rm (2)}] $\pp^2\dasharrow \pp^2$ with defining coordinates $(xz:yz:y^2)$

The base locus of this Cremona map consists of the points $(0:0:1)$ and $(1:0:0)$,
with multiplicity $1$ and $2$, respectively  (in the classical terminology, one proper point and another
proper point with a point in its first neighborhood).

\item[{\rm (3)}] $\pp^2\dasharrow \pp^2$ with defining coordinates $(y^2-xz:yz:z^2)$

The base locus of this Cremona map consists of the point $(1:0:0)$ with multiplicity $3$,
a so-called triple structure on a point
(in the classical terminology, one proper point with a point in its first neighborhood
and a point in its second neighborhood).
\end{enumerate}

\begin{Theorem}\label{rational quartics}
{\rm ($k$ algebraically closed of characteristic zero)}
Let $F=F(\uu,x,y,z)\in k[\uu,x,y,z]$ be a family of rational
plane curves of degree $4$ with a fixed set of singular points in the sense previously defined.
Then the general member $f=F(\alpha)\in k[x,y,z]$ in this family has
gradient ideal of linear type.
\end{Theorem}
\demo
We will actually show a bit more, namely, that any irreducible rational quartic falls within
a family whose general member has the required property for its gradient ideal.
In this vein, we can and will assume that the members of any family are singular.
This is because the gradient ideal of any smooth plane curve is generated by a regular sequence,
hence trivially of linear type.

Now, any irreducible rational quartic $f\in k[x,y,z]$ has at least one double singular point and at most a triple point.
Let us first consider the situation where $f$ has a double point -- hence has at most $3$ such
points and no triple point.

In this case, as explained above, $f$ comes from a conic by means of a Cremona map.

Let $Q=u_1x^2+u_2y^2+u_3z^2+2u_4yz+2u_5zx+2u_6xy$
be the equation of the conic as above, assumed non-singular, i.e.,
the corresponding symmetric matrix has nonzero determinant $\Delta=u_1u_2u_3+2u_4u_5u_6-u_1u_4^2-u_2u_5^2-u_3u_6^2$
(the {\em discriminant} of $Q$).

Applying the above Cremona maps, we obtain, respectively:
\begin{enumerate}
\item[{\rm (1)}]
A quartic with exactly three double points at $P_1=(1:0:0)$, $P_2=(0:1:0)$ and $P_3=(0:0:1)$,
where $P_1$ (respectively, $P_2$, $P_3$) is a node except when
the principal minor $u_1u_2-u_6^2$ vanishes (respectively, except when the principal minors
$u_1u_3-u_5^2$, $u_2u_3-u_4^2$ vanish).

Here we may harmlessly assume that $u_1=u_2=u_3=1$
provided they are all nonzero.

In this block belong the following families.
\begin{enumerate}
\item[{\rm (a)}] Three nodes:
$$\tilde{f}= y^2z^2+x^2z^2+x^2y^2+2xyz(u_4x+u_5y+u_6z),\ \ \ \ u_4,u_5,u_6\neq\pm 1, \ \  \Delta\neq 0$$
where $\Delta=2u_4u_5u_6-u_4^2-u_5^2-u_6^2+1.$
\end{enumerate}
\begin{enumerate}
\item[{\rm (b)}] Two nodes and one cusp:
$$\tilde{f}= y^2z^2+x^2z^2+x^2y^2+2xyz^2+2xyz(u_4x+u_5y),\ \ \ \  u_4,u_5\neq\pm 1,\ \   \Delta\neq 0$$
where $\Delta=2u_4u_5-u_4^2-u_5^2=-(u_4-u_5)^2$.
\end{enumerate}
\begin{enumerate}
\item[{\rm (c)}] One node and two cusps:
$$\tilde{f}= y^2z^2+x^2z^2+x^2y^2+2xy^2z+2xyz^2+2u_4x^2yz,  \ \ \ u_4\neq \pm 1$$
\end{enumerate}
\begin{enumerate}
\item[{\rm (d)}] Three cusps :
$$\tilde{f}= y^2z^2+x^2z^2+x^2y^2-2xyz(x+y+z)$$
\end{enumerate}

\item[{\rm (2)}]  A quartic with a double point at $P_1=(0:0:1)$ (which is a node or a cusp according
as to whether $u_2\neq 0$ or $u_2=0$) and a double point at $P_2=(1:0:0)$
(which is either a tacnode or a ramphoid cusp according as to whether the principal minor $u_1u_3-u_5^2$
is nonzero or vanishes).

Here we may assume that $u_1=u_3=1$ and $u_6=0$.

It comprises the following families.
\begin{enumerate}
\item[{\rm (e)}] One tacnode and one cusp:
$$\tilde{f}= x^2z^2+y^4+2y^3z+2u_5xy^2z, \ \ \ \ \ \ \ \ \ \  \ \ \ u_5\neq \pm 1$$
\end{enumerate}
\begin{enumerate}
\item[{\rm (f)}] One tacnode and one node :
$$\tilde{f}= z^2(x^2+y^2)+y^4+2y^2z(u_4y+2u_5x),  \ \ \ \ u_5\neq \pm 1,\ \ \ \Delta=u_4^2+u_5^2-1\neq 0$$
\end{enumerate}
\begin{enumerate}
\item[{\rm (g)}] One ramphoid cusp  and one node:
$$\tilde{f}= x^2z^2+y^4+2zy^3+2xy^2z+u_2z^2y^2, \ \ \ \ u_2\neq 0$$
\end{enumerate}
\begin{enumerate}
\item[{\rm (h)}] One ramphoid cusp  and one cusp:
$$\tilde{f}= x^2z^2+y^4+2zy^3+2xy^2z$$
\end{enumerate}

\item[{\rm (3)}] A quartic with  an oscnode at $P_1=(1:0:0)$ if $u_2\neq 0$; else, a
singularity of type $A_6$.

Here we may assume that  $u_1=1$ and $u_5=u_6=0$.
Namely, we get the following forms.

\begin{enumerate}
\item[{\rm (i)}] One oscnode :
$$\tilde{f}= (y^2-xz)^2+y^2z^2+u_3z^4,\ \ \ \ \ \ \ \ \ \  u_3\neq 0$$
\end{enumerate}
\begin{enumerate}
\item[{\rm (j)}] One singularity of type $A_6$:
$$\tilde{f}= (y^2-xz)^2+2yz^3$$
\end{enumerate}

\end{enumerate}
An irreducible rational quartic having only double points -- ordinary or not -- falls
within the following families up to coordinate change, according to the nature of its singularities.
We have written $\tilde{f}$ instead of $F$ to help us think of the general member instead of the family
itself.
To keep track of the parameters in each case we have maintained the original indices, however anaesthetical
they may look.

Note the two kinds of degeneration:  first, whether a member of the family factors is controlled
by the vanishing of the corresponding value of $\Delta$; second, whether the non-general member goes across
stratified  subfamilies is controlled by the vanishing of another ideal in $k[\uu]$ - we will
call the latter ideal the {\em strata degeneration locus}.

\medskip

We now consider the case where the quartic has a triple point, say, at $P=(0:0:1)$.

Generally, for a plane irreducible curve of degree $d$ with a singular point of multiplicity
$d-1$ (hence, a rational curve), it is frequently easier to look at the linear type condition.
In the case of a quartic, up to a projective change of coordinates,  the
equation of the curve has the form $\phi(x,y)\, z+\psi(x,y)=0$,
where $\phi$ can moreover be brought up to one of the forms  $x(y^2-x^2)$, $xy^2$ and $y^3$, and $\psi$
may be further normalized in such a way that the resulting family has as few parameters
as possible.

\medskip

(4)
After these reductions,  any irreducible
plane quartic having $(0:0:1)$ as a triple point falls within three basic families, according to the nature
of the triple point:
\begin{enumerate}
\item[{\rm (k)}] An ordinary triple point:
$$\tilde{f}= x(y^2-x^2)z+y^4+x^2y(u_1y+u_2x) $$
\end{enumerate}
\begin{enumerate}
\item[{\rm (l)}] A triple point with double tangent:
$$\tilde{f}= xy^2z+x^4+y^4+u_1x^3y$$
\end{enumerate}
\begin{enumerate}
\item[{\rm (m)}] A higher cusp :
$$\tilde{f}= y^3z+x^4+u_1x^2y^2.$$
\end{enumerate}

\bigskip

The proof proceeds by dealing with each of the above types.

\bigskip

{\sc Block (1)}

\medskip

(a) By Theorem~\ref{degeneration_vs_lineartype}, it suffices to show that $\mathcal{I}({\bf 0})$ has codimension $3$.
From the symmetrical parametric structure of $\tilde{f}$, the inclusion $\mathcal{I}({\bf 0})\subset I_{\tilde{f}({\bf 0})}$
is an equality.
On the other hand, by direct verification, the latter ideal admits the following syzygies:
$$\left(
  \begin{array}{c}
    xy^2-xz^2\\
    -y^3-yz^2\\
    y^2z+z^3
  \end{array}
\right),\quad
\left(
  \begin{array}{c}
  -x^3-xz^2\\
  x^2y-yz^2\\
  x^2z+z^3
 \end{array}
\right)
$$
Bringing in besides  the generators of the gradient ideal, one obtains after a calculation
 $x^3,y^3,z^3 \in\mathcal{I}({\bf 0})$.
Thus, $\mathcal{I}({\bf 0})$ has codimension $3$.

\medskip

We first note that (b)--(d) are obvious successive strata of the family (a).

(b) In this first stratum  the evaluated ideal $\mathcal{I}({\bf 0})$ has codimension $2$,
so one may try another evaluation.
Note that one cannot blindly apply the result of (a) to claim that, here too, there is some general value
$\alpha$ for which the assertion holds, since $\alpha$ could lie outside the open set obtained in (a).
(Of course, it goes without saying that an explicit such open set can be computed by way of
Theorem~\ref{degeneration_vs_lineartype}, (ii) and we even have a conjecture about its form.)
Instead, we resort to a painful hand verification, namely, the following is a syzygy of
$\mathcal{I}$:
$$
\left(
  \begin{array}{c}
    x^2(u_4^2-1)+xy(u_4u_5-1)+2xz(u_4-u_5)+yz(u_4-u_5) \\
    y^2(u_5^2-1)+ xy(u_4u_5-1)+xz(u_5-u_4)+2yz(u_5-u_4)\\
    3z^2(u_4-u_5)+ xy(u_4-u_5)+xz(2u_4^2-u_5u_4-1)+yz(-2u_5^2+u_5u_4+1) \\
  \end{array}
\right)
$$
Looking at the first summand of each coordinate, one sees that for
any ``value'' $\alpha=(u_4,u_5)\in \mathbb{A}_k^2$
with $u_4\neq u_5$ and $u_4\neq\pm1$, $u_5\neq\pm 1$,
pure powers of $x,y,z$ remain and the ideal generated by the coordinates will
have codimension $3$, hence also the corresponding  $\mathcal{I}(\alpha)$.

\medskip

(c) In this stratum by the same token, we look at the following two syzygies:

$$\left(
  \begin{array}{c}
 -xy+xz\\
3y^2+2xyu_4+xz+3yz\\
-3z2-2xzu_4-2xy-3yz
  \end{array}
\right),\quad
\left(
  \begin{array}{c}
x^2(u_4+1)+3/2xy+3/2xz+yz \\
    3/2y^2-xy(u_4+1)+1/2yz \\
    3/2z^2+2/2yz+xz(u_4+1)
 \end{array}
\right)
$$
An identical analysis as above, looking at the pure powers, allow to choose any
``value'' $u_4\neq -1$.

\medskip

(d)
Here too it suffices to look at the following syzygies
$$\left(
  \begin{array}{c}
  x^2+xy+2/3xz-2/3yz\\
-xy-y^2+2/3xz-2/3yz\\
-1/3xz+1/3yz
\end{array}
\right),\quad
\left(
  \begin{array}{c}
xy+xz-2/3yz\\
-y^2+1/3yz\\
1/3yz-z^2
 \end{array}
\right)
$$
Once more, pure powers of the variables are easily located.

\bigskip

{\sc Block (2)}

\medskip

(e)
As in the proof of (a), here too the ideal  $\mathcal{I}({\bf 0})$ contains the ideal of $1$-minors of the syzygy matrix
of $I_{\tilde{f}({\bf 0})}$.
One can check that the latter ideal admits the following syzygies:
$$\left(
  \begin{array}{c}
   2x^2-3y^2\\
xz\\
-2xz
  \end{array}
\right),\quad
\left(
  \begin{array}{c}
  2xy+3xz\\
 yz\\
 -2yz-3z2
 \end{array}
\right).
$$

Bringing in the generator $\partial \tilde{f}({\bf 0})/\partial y = 4y^3+6y^2z$, one readily
sees that $\mathcal{I}({\bf 0})$ has codimension $3$.

\medskip

(f) This follows the same pattern as (e).
The relevant syzygies to look at are
$$\left(
  \begin{array}{c}
 x^2+y^2\\
0\\
-xz
  \end{array}
\right),\quad
\left(
  \begin{array}{c}
 2xy^2+xz^2\\
yz^2\\
-2y^2z-z3
 \end{array}
\right)
$$
and the calculation to get suitable powers of $x,y,z$ inside $\mathcal{I}({\bf 0})$
is pretty straightforward.

\medskip

(g)
If we indulge ourselves allowing a computation with {\em Macaulay}, we
get $u_2^2\in (\mathcal{I}\colon \mathfrak{m}^{\infty})$.
By Theorem~\ref{degeneration_vs_lineartype}, {\em every} member in this family has gradient ideal of linear type.
Alternatively, one can look for $\mathcal{I}$ evaluated, say, at $u_2\mapsto 1$.
One column turns out to be
$$
\left(
  \begin{array}{c}
5x^2-y^2+xz-yz\\
y^2+xy+xz+yz\\
-z^2-2y^2-xz-2yz
 \end{array}
\right)
$$

Thus, the ideal generated by the three coordinates above already has codimension $3$.

\medskip

(h)
Again, a computation with {\em Macaulay} gives the following syzygies:
$$\left(
  \begin{array}{c}
   2y^2-3xz\\
-yz\\
3z^2
  \end{array}
\right),\quad
\left(
  \begin{array}{c}
  x^2-27/50xz\\
1/5xy+1/5y^2+3/25xz-9/50yz\\
-2/5y^2-xz-6/25yz+27/50z^2
 \end{array}
\right).
$$
Thus, we locate pure powers as terms of the coordinates.
Alternatively, note the coordinates $3z^2$ and $x^2-27/50xz=x(x-27/50z)$ are invertible locally at $(x,y)$ and
$(y,z)$, respectively. Since the latter are the two singular minimal primes of the quartic,
this shows that the gradient ideal is locally a complete intersection
at these primes, hence is of linear type.

\bigskip

{\sc Block (3)}

\medskip

(i) The discussion of this case is analogous to the one of (g). A computation with {\em Macaulay}  yields
get $u_3\in (\mathcal{I}\colon \mathfrak{m}^{\infty})$.
Therefore, {\em every} member in this family has gradient ideal of linear type.
As previously enacted,  the ideal $\mathcal{I}(1)$
has codimension $3$ through a convenient analysis of its terms.

\medskip

(j)
The following vectors are directly seen to be syzygies of the gradient ideal:
$$
\left(
  \begin{array}{c}
  6y^2+xz \\
3yz    \\
-z^2
\end{array}
\right),\quad
\left(
  \begin{array}{c}
 7x^2+18yz \\
3xy\\
6y^2-7xz
\end{array}
\right).
$$
A straightforward calculation gives the right codimension.

\bigskip

{\sc Block (4) (Triple point)}

\medskip

(k) The discussion of this case is much like the one of (a)
in that the parametric structure of $\tilde{f}$
allows to see that the syzygies of $I_{\tilde{f}({\bf 0})}$ are contained
in the syzygies of $I_{\tilde{f}}$ evaluated at ${\bf 0}$.
It then suffices to check that $I_{\tilde{f}({\bf 0})}$ is of linear type.
We do this by arguing that it is locally a complete intersection at its unique singular
prime $(x,y)$.
For this, it suffices to consider the syzygy
$$\left(
  \begin{array}{c}
x^2-2/3y^2+1/6xz\\
1/6yz\\
-(3x+1/2z)z
\end{array}
\right)
$$
where the last coordinate is invertible locally at $(x,y)$.

\medskip

(l) As in (k), the syzygies of $I_{\tilde{f}({\bf 0})}$ are contained
in the syzygies of $I_{\tilde{f}}$ evaluated at ${\bf 0}$.
Here it is elementary to guess the syzygy
$$\left(
  \begin{array}{c}
0\\
xy\\
-4y^2-2xz
\end{array}
\right)
$$
Using further the generators of the gradient ideal of $\tilde{f}({\bf 0})$, it is readily
seen powers of the variables among the entries.

\medskip

(m) This case is like the previous one, only more elementary.
We argue that $\mathcal{I}({\bf 0})$ has codimension $3$ as before
by looking at the obvious syzygy of $I_{\tilde{f}({\bf 0})}$
$$\left(
  \begin{array}{c}
0\\
y\\
-3z
\end{array}
\right)
$$
which clearly tells us that the ideal is locally a complete
intersection at $(x,y)$.
\qed

\begin{Remark}\rm
We have drawn quite a bit on computation to verify all cases of the theorem.
Using Theorem~\ref{degeneration_vs_lineartype}, to have a computation-free argument it would suffice to show that
the contraction $(\mathcal{I}\colon \mathfrak{m}^{\infty})\cap k[\uu]$
coincides set-theoretically with the product $\Delta\mathfrak{a}$ of the discriminant
and the strata degeneration locus $\mathfrak{a}\subset k[\uu]$.
We conducted a computational verification of this fact for all four blocks
of families. Thus, morally, the conjecture for rational quartics is settled.
However, we have found no immediate theoretical reason pointing at least to an ideal inclusion
$(\mathcal{I}\colon \mathfrak{m}^{\infty})\cap k[\uu]\subset \Delta$.
Note that, according to Proposition~\ref{when_cod_is_one}, one expects that
the contraction also have codimension one, whereas the strata degeneration locus
is given by a principal ideal.

It seems natural to conjecture that any {\em irreducible} quartic has gradient ideal of linear type.
\end{Remark}

{\bf Authors addresses:}

\medskip

\noindent {\sc Departamento de Matem\'atica, CCEN, Universidade Federal
de Pernambuco,
Cidade Universit\'aria, 50740-540 Recife, PE, Brazil}

\noindent E-mail:  abbasnn@dmat.ufpe.br, aron@dmat.ufpe.br

%\bla

\end{document}